\documentstyle[12pt,leqno,amssymb]{article}

\newtheorem{prop}{Proposition}[section]

\def\R{{\bf R}}
\def\C{{\bf C}}
\def\Z{{\bf Z}}
\def\Q{{\bf Q}}

\def\A{{\bf A}}

\def\a{ \alpha }
\def\b{ \beta }
\def\ba{ \backslash}
\def\c{ \chi }
\def\d{ \delta }
\def\e{ \eta }
\def\g{ \gamma}

\def\D{ \Delta }
\def\l{ \lambda }
\def\L{ \Lambda }

\def\s{ \sigma }

\def\p{ \varphi }

\def\z{ \zeta }

\def\no{ \noindent }

\def\w{ \varpi }
\def\ra{ \rightarrow }
\def\Re{{\rm Re}}
\def\Im{{\rm Im}}
\def\Pic{{\rm Pic}}
\def\Hom{{\rm Hom}} 
\def\Ker{{\rm Ker}} 
\def\Val{{\rm Val}} 
\def\Spec{{\rm Spec}}
\def\frak{\bf } 
\def\Norm{{\rm Norm}}
\def\Vol{{\rm Vol}}
\def\tp{\tilde{\varphi}}
\def\tl{\tilde{\lambda}}
\def\U{{\cal U}}
\def\SL{{\Sigma_1- \Sigma_1'}}
\def\DL{{\Delta_P -\Delta_{P'}}}
\def\ce{{\chi}_{\eta}}
\def\Ind(L){{\cal I}}
\def\ind(L){{\iota}}
\def\KX{{\omega}^{\vee}}

\title{Height zeta functions
of \\
 toric bundles over flag varieties}

\author{Matthias Strauch\\
\small Mathematisches Institut, Universit\"at Bonn \\
\small Beringstr. 1, 53117 Bonn, Germany\\
\small e-mail: mats@rhein.iam.uni-bonn.de\\
and\\
Yuri Tschinkel\\
\small Dept. of Mathematics, U.I.C.\\
\small Chicago, (IL) 60608,  U.S.A.  \\
\small e-mail: yuri@math.uic.edu
}

\begin{document} 
 
\date{}

\maketitle

\thispagestyle{empty}

\pagebreak
\tableofcontents

\bigskip
\bigskip
\bigskip

\section{Introduction}
\label{1}

\subsection*
\no
{\bf 1.1}\hskip 0,5cm
Let $X$ be a nonsingular projective algebraic variety over a number field
$F$. Let ${\cal L}=\left(L, (\|\cdot\|_v)_v\right)$ be a  metrized line 
bundle on $X$, i.e., a line bundle $L$ together with a family of 
$v$-adic metrics, where $v$ runs over the set $\Val (F) $ of places of $F$.
Associated to ${\cal L}$ there is a height function
$$
H_{{\cal L}}\,:\, X(F)\ra \R_{>0}
$$
on the set $X(F)$ of $F$-rational points of $X$ (cf. \cite{ST,Peyre}
for the definitions of $v$-adic metric, metrized line bundle and height
function). 
For appropriate subvarieties $U\subset X$ and line bundles $L$ 
we have
$$
N_U({\cal L},H):=\{x\in U(F)\,\,|\,\,H_{\cal L}(x)\le H\} <\infty
$$
for all $H$ (e.g., this holds for any $U$ if $L$ is ample).
We are interested in the asymptotic behavior of this counting function 
as $H\ra \infty$. It is expected that the behavior of such 
asymptotics can be described in geometric terms (\cite{BaMa,FMT}).

\no
Let 
$$
\L_{\rm eff}(X):=\sum_{H^0(X,L)\neq 0} \R_{\ge 0}[L]\subset \Pic (X)_{\R}
$$
be the closed cone in $\Pic(X)_{\R}$ generated by the classes of effective
divisors ($[L]$ denotes the class of the line bundle $L$ in $\Pic(X)$).
Let $L$ be a line bundle on $X$ such that $[L]$ lies in the interior of 
$\L_{\rm eff}(X)$. Define 
$$
a(L):=\inf\{a\in \R\,\,|\,\, a[L]+[K_X]\in \L_{\rm eff}(X)\},
$$
where $K_X$ denotes the canonical line bundle on $X$.
Assume that $\L_{\rm eff}(X)$ is a finitely generated polyhedral cone.
For $L$ as above we let $b(L)$ be the codimension of the minimal face
of $\L_{\rm eff}(X)$ which contains  $a(L)[L]+[K_X]$.

By a Tauberian theorem 
(cf. \cite{De}, Th\'eor\`eme III),
the asymptotic behavior of $ N_U({\cal L},H)$ can be determined if one has
enough information about the height zeta function 
$$
Z_U({\cal L},s):= \sum_{x\in U(F)}H_{\cal L}(x)^{-s}.
$$
More precisely, 
suppose that $Z_U({\cal L},s)$ converges for $\Re(s)\gg 0$,
that it  has an abscissa
of convergence $a>0$ and that it 
can be continued meromorphically to a half-space
beyond the abscissa of convergence. 
Suppose further that there is a pole of 
order $b$ at $s=a$ and that there are no other poles in this half-space.
Then
$$
N_U({\cal L},H)=cH^a(\log H)^{b-1}(1+o(1))
$$
for $H\ra \infty$ and 
$$
c=\frac{1}{(b-1)!a}\cdot \lim_{s\ra a} (s-a)^bZ_U({\cal L},s).
$$
It is conjectured that for appropriate $U$ and ${\cal L}$ 
one has $a=a(L)$ and $b=b(L)$ (cf.  \cite{BaMa,FMT}).
Moreover, there is a conjectural framework how to 
interpret the constant $c$ (cf. \cite{Peyre,BaTschi5}). 
There are examples which demonstrate that
this geometric ``prediction'' of the asymptotic 
cannot hold in complete generality, even for smooth
Fano varieties  (cf. \cite{BaTschi4}).
Our goal is to show that
the conjectures  do hold for a class of varieties closely related to
linear algebraic groups. Our results are a natural extension of
corresponding results for flag  varieties (cf. \cite{FMT}) and 
toric varieties (cf. \cite{BaTschi2,BaTschi3}).  We proceed to describe 
the class of varieties under consideration.

\subsection*
\no
{\bf 1.2}\hskip 0,5cm
Let $G$ be a semi-simple simply connected split algebraic group
over $F$ and $P\subset G$ an $F$-rational parabolic subgroup of $G$.
Let $T$ be a split algebraic torus over $F$ and $X$ a projective 
nonsingular equivariant compactification of $T$.
A homomorphism $\eta\,:\, P\ra T$ gives rise to an action of $P$
on $X\times G$ and the quotient $Y:=(X\times G)/P$ is again a nonsingular
projective variety over $F$. There is a canonical morphism 
$\pi\,:\, Y\ra W:= P\ba G$ such that $Y$ becomes a locally trivial fiber 
bundle over $W$ with fiber $X$. 
Corresponding to a character $\l\in X^*(P)$ there is a line bundle $L_{\l}$
on $W$ and the assignment $\l\mapsto L_{\l}$ gives an isomorphism 
$X^*(P)\ra \Pic(W)$. 

The toric variety 
$X$ can be described combinatorially by a fan $\Sigma$ in the
dual space of the space of characters $X^*(T)_{\R}$. 
Let $PL(\Sigma)$ be the group 
of $\Sigma$-piecewise linear integral 
functions on the dual space of $X^*(T)_{\R}$.
Any $\p\in PL(\Sigma)$ defines a line bundle 
$L_{\p}$ on $X$ which is equipped 
with a canonical $T$-linearization and we get an isomorphism 
$PL(\Sigma)\simeq \Pic^T(X). $ There is a canonical exact sequence
$$
0\ra X^*(T)\ra PL(\Sigma)\ra \Pic(X)\ra 0.
$$
The $T$-linearization of $L_{\p}$ allows us to define a line bundle
$L^Y_{\p}$ on $Y$ and this gives a homomorphism $PL(\Sigma)\ra \Pic(Y)$.
One can show that there is an exact sequence
\begin{equation}
\label{1.2.1}
0\ra X^*(T)\ra PL(\Sigma)\oplus X^*(P)\ra \Pic(Y)\ra 0.
\end{equation}
Denote by $Y^{o}:=(T\times G)/P$ the open subvariety of $Y$ obtained as
the twist of $T $ with $W$.

\subsection*
\no
{\bf 1.3} \hskip 0,5cm
By means of a maximal compact subgroup in
the adelic group $G({\bf A})$ we can introduce
metrics on the line bundles $L_{\l}$. The corresponding height zeta 
functions are Eisenstein series:
$$
\sum_{w\in W(F)}H_{{\cal L}_{\l}}(w)^{-s}=E^G_P(s\l-\rho_P,1_G).
$$
On the other hand, for any $\p\in PL(\Sigma)$ there is a function
$$
H_{\Sigma}(\,\cdot\, ,\p)\,:\, T({\bf A})\ra \R_{>0}
$$
such that $H_{\Sigma}(x,\p)^{-1}$ is the height of $x\in T(F)$ 
with respect to a metrization ${\cal L}_{\p}$ of $L_{\p}$.
This metrization induces a metrization ${\cal L}_{\p}^Y$ of
the line bundle $L^Y_{\p}$ on $Y$.

Let $(x,\g)\in T(F)\times G(F)$  and let $y$ be the image of 
$(x,\g) $ in $Y(F)$. Then there is a $p_{\g}\in P({\bf A})$ such 
that
$$
H_{{\cal L}^Y_{\p}}(y)=H_{\Sigma}(x\eta(p_{\g}),\p)^{-1}.
$$
Hence we may write formally
\begin{equation}
\label{1.3.3}
Z_{Y^o}({\cal L}_{\p}^Y\otimes \pi^*{\cal L}_{\l},s)
=\sum_{\g\in P(F)\ba G(F)}H_{{\cal L}_{\l}}(\g)^{-s}\sum_{x\in T(F)}
H_{\Sigma}(x\eta(p_{\g}),s\p).
\end{equation}
Now we apply Poisson's summation formula for  the 
the torus and get
\begin{equation}
\label{1.3.4}
\sum_{x\in T(F)}
H_{\Sigma}(x\eta(p_{\g}),s\p)=\int_{(T({\bf A})/T(F))^*}\hat{H}_{\Sigma}(\chi,s\p)
\chi(\eta(p_{\g}))^{-1}d\c,
\end{equation}
where $\hat{H}_{\Sigma}(\,\cdot\, ,s\p)$ denotes the Fourier transform
of ${H}_{\Sigma}(\,\cdot\, ,s\p)$ and  $(T({\bf A})/T(F))^*$ is 
the group of unitary characters of 
$T({\bf A})$ which are trivial on $T(F)$
equipped with the orthogonal measure 
$d\c$. Actually, it is sufficient to consider 
only those characters which are trivial on the maximal compact subgroup
${\bf K}_T$ of $T({\bf A})$, because the function $H_{\Sigma}(\,\cdot\, , \p)$ is
invariant under ${\bf K}_T$.
The expression (\ref{1.3.4}) can now be put into (\ref{1.3.3}), and 
after interchanging summation and integration the result is
\begin{equation}
\label{1.3.5}
Z_{Y^o}({\cal L}_{\p}^Y\otimes \pi^* {\cal L}_{\l},s)=
\int_{(T({\bf A})/T(F){\bf K}_T)^*}\hat{H}_{\Sigma}(\c,s\p)E^G_P(s\l-\rho_P,
(\chi\circ\eta)^{-1})d\c
\end{equation}
where $E^G_P(s\l-\rho_P,\xi)=
E^G_P(s\l-\rho_P,\xi,1_G)$ is the Eisenstein series
twisted by a character $\xi$ of $ P({\bf A})$. 
This is the starting point for the investigation of the height zeta function. 
To get an expression which is more suited  for our study we decompose
the group of characters $(T({\bf A})/T(F){\bf K}_T)^*$
into a continuous and a discrete part, i.e.,
$$
(T({\bf A})/T(F){\bf K}_T)^*=X^*(T)_{\R}\oplus {\cal U}_T,
$$
where $X^*(T)_{\R}$ is the continuous part and ${\cal U}_T$
is the discrete part. 
The right-hand side of (\ref{1.3.5}) is accordingly
\begin{equation}
\label{1.3.6}
\int_{X^*(T)_{\R}}\left\{ \sum_{\c\in {\cal U}_T}
\hat{H}_{\Sigma}(\c,s\p+ix)E^G_P(s\l-\rho_P -i\check{\eta}(x),
(\chi\circ\eta)^{-1})\right\}dx.
\end{equation}
Recall that we would like to show that this function which is
defined for $\Re(s)\gg 0$ (assuming that $(\p,\l)$ is
contained in a convex open cone) can be continued meromorphically
beyond the abscissa of convergence. To achieve this we need more
information on the function under the integral sign in (\ref{1.3.6}).
First we have to determine the singularities of 
$$
(\p,\l)\mapsto \hat{H}_{\Sigma}(\c,\p)E^G_P(\l-\rho_P,(\c\circ\eta)^{-1})
$$
near the cone of absolute convergence.
This is possible because $\hat{H}_{\Sigma}$
can be calculated rather explicitly and it is not so difficult to 
determine the singular hyperplanes of the Eisenstein series with characters.
The next step consists in an iterated application of Cauchy's residue 
formula to the integral over the real vector space $X^*(T)_{\R}$.
This can be done only if one knows that
\begin{equation}
\label{1.3.7}
\sum_{\c\in {\cal U}_T}
\hat{H}_{\Sigma}(\c,s\p+ix)E^G_P(s\l-\rho_P-i\check{\eta}(x),
(\c\circ\eta)^{-1})
\end{equation}
satisfies some growth conditions when $x\in X^*(T)_{\R}$ tends to 
infinity. This is true for the function 
$x\mapsto \hat{H}_{\Sigma}(\c,s\p+ix)$ thanks to the explicit expression
mentioned above. The absolute value of the Eisenstein series 
$E^G_P(s\l-\rho_P -i\check{\eta}(x),(\c\circ\eta)^{-1})$ 
will in general increase for $x\ra \infty$
if $\Re(s)\l-\rho_P$ is not contained in the cone of absolute 
convergence. However, if $\Re(s)\l-\rho_P$ 
is sufficiently close to the boundary
of that cone, this increasing behavior is absorbed by the decreasing 
behavior of $\hat{H}_{\Sigma}(\c,s\p+ix)$.

Therefore, we may apply Cauchy's residue theorem and show that
(\ref{1.3.5}) can be continued meromorphically
to a larger half-space and that there are no poles (in $s$) with
non-zero imaginary part.

The Tauberian Theorem can now be used to prove asymptotic
formulas for the counting function 
$N_{Y^o}({\cal L}^Y_{\p}\otimes \pi^*{\cal L}_{\l},H)$ provided that one 
knows the order of the pole of the height zeta function.
This problem can be reduced to the question whether the ``leading term''
of the Laurent series of (\ref{1.3.7}) does not vanish.
That this is indeed so will be shown in section 6.

\subsection*
\no
{\bf 1.4}\hskip 0,5cm
We have restricted ourself to the case of split tori and split groups
because this simplifies some technical details. The general case
can be treated similarly. 
 
We consider these results as an important step towards an understanding
of the arithmetic of spherical varieties. For example, choosing $P=B$ a Borel 
subgroup, $T=B/U$ where $U$ is the unipotent radical of $B$ and 
$\eta\,:\, B\ra T$ the natural projection, we obtain an equivariant 
compactification of $U\ba G$, a horospherical variety.

We close this introduction with a brief description of the remaining
sections. 
Section 2  recalls the relevant facts we need concerning generalized flag
varieties, i.e., description of line bundles on $W=P\ba G$, the
cone of effective divisors in $\Pic(W)_{\R}$, metrization of line bundles,
height zeta functions. The exposition is based entirely on the paper 
\cite{FMT}.

The next section contains the corresponding facts for toric varieties. It
is a summary of a part of \cite{BaTschi1}.
We give the explicit calculation of the Fourier transform 
$\hat{H}_{\Sigma}(\,\cdot \, ,\p)$ and show 
that Poisson's summation formula can 
be used to give an expression of the height zeta function 
$Z_T({\cal L}_{\p},s)$ .

In section 4 we introduce twisted products, discuss line bundles on these,
the Picard group (cf. (\ref{1.2.1})), metrizations of line bundles etc.
It ends with the formula (\ref{1.3.5}) for the height zeta function 
$Z_{Y^o}({\cal L}_{\p}^Y\otimes \pi^* {\cal L}_{\l},s)$
in the domain of absolute convergence. 

The first part of section 5 explains the method for the proof that the
height zeta function can be continued meromorphically to a half-space
beyond the abscissa of absolute convergence. Moreover, we state a theorem
which gives a description of the coefficient of the Laurent series  at the pole
in question. This coefficient will be the leading one, provided that it does
not vanish. One can relate the coefficient to arithmetic and
geometric invariants
of the pair  $(U,{\cal L})$ but we decided not to pursue this,
since there are detailed expositions of all the necessary arguments
in \cite{Peyre,BaTschi1,BaTschi5}.

These two theorems (meromorphic continuation of certain
integrals and the description of the coefficient) 
will be proved in a more general context in section 7.
The second part of section 5 contains the proof that the hypothesis of these
theorems are fulfilled in our case.
It ends with the main theorem on the asymptotic behavior of the counting
function $N_{Y^o}({\cal L},H)$, assuming that the coefficient of the
Laurent series mentioned above does not vanish. 
Section 6 is devoted to the proof of this fact.
In the last section we prove some statements on Eisenstein series (well-known
to the experts) which are used in section 5.

\bigskip
\no
{\bf Acknowledgements.}
We are very grateful to V. Batyrev and J. Franke for helpful
discussions and collaboration on related questions. 
The first author was supported by the DFG-Graduiertenkolleg of the
Mathematics Institute of the University of Bonn. 
Part of this work was done while the second author was visiting
the MPI in Bonn, ETH Z\"urich and ENS Paris. He would like to thank
these institutions for their hospitality.

\bigskip
\no
{\bf Some notations.} \hskip 0,5cm
In this paper $F$ always denotes a fixed algebraic number field.
The set of places of $F$ will be denoted by $\Val(F)$ and the
subset of archimedean places by $\Val_{\infty}(F)$. We shall 
write $v\mid\infty$  if $v\in \Val_{\infty}(F)$ and
$v\nmid \infty$ if $v\notin \Val_{\infty}(F)$ . 
For any place $v$ of $F$ we denote by $F_v$ 
the completion of $F$ at $v$ and by ${\cal O}_v$ the ring of 
$v$-adic integers (for $v\nmid \infty $).  
The local absolute value $|\cdot|_v$ on $F_v$ is
the multiplier of the Haar measure, i.e., $d(ax_v)=|a|_vdx_v$ for some
Haar measure $dx_v$ on $ F_v$.
Let $q_v$ be the cardinality of the residue field of $F_v$
for non-archimedean valuations
and put $q_v=e$ for archimedean valuations. We denote by ${\bf A}$
the adele ring of $F$. 
 For any algebraic group
$G$ over $ F$ we denote by $X^*(G)$ the group of (algebraic) characters
which are defined over $F$.

\bigskip
\bigskip
\bigskip

\section{Generalized flag varieties}
\label{2}

\subsection*
\no
{\bf 2.1}\hskip 0,5cm
Let $G$ be a semi-simple simply connected linear algebraic group
which is defined and split over $F$. We fix a Borel subgroup
$P_0$ over $F$ and a Levi decomposition $P_0=S_0U_0$ with a maximal 
$F$-rational 
torus $S_0$ of $G$. Denote by $P$ a standard (i.e., containing $P_0$) 
parabolic subgroup and by $W=P\backslash G$ the corresponding flag
variety. The quotient morphism $G\ra W$ will be denoted by $\pi_W$. 
Any character $\l\in X^*(P)$ defines a line bundle $L_{\l}$
on $W$ by
$$
\Gamma(U,L_{\l}):=
\{f\in {\cal O}_G(\pi^{-1}_W(U))\,\,|\,\, f(pg)=\l(p)^{-1}f(g)
\hskip 0,2cm \forall g\in \pi^{-1}_W(U), p\in P\, \}.
$$
The assignment $\l\mapsto L_{\l}$ gives an isomorphism (because $G$ is assumed
to be simply connected)
$$
X^*(P)\ra {\rm Pic}(W)
$$
(cf. \cite{Sa}, Prop. 6.10). The anti-canonical line bundle $\KX_W$
corresponds to $2\rho_P$ (the sum of roots of $S_0$ 
occurring in the unipotent radical of $P$.)

\subsection*
\no 
{\bf 2.2}\hskip 0.5cm
These line bundles will be metrized as follows. Choose a maximal compact 
subgroup ${\bf K}_G=\prod_{v}{\bf K}_{G,v}\subset G({\bf A})$ 
(${\bf K}_{G,v}\subset G(F_v)$), such that the Iwasawa decomposition
$$
G({\bf A})=P_0({\bf A}){\bf K}_G
$$
holds.
Let $v\in {\rm Val}(F)$ and $w\in W(F_v)$. 
Choose $k\in {\bf K}_{G,v}$ which is 
mapped to $w$ by $\pi_W$. For any local 
section $s$ of $L_{\l}$ at $w$ 
we define
$$
\|w^*s\|_w:=|s(k)|_v.
$$
This gives a $v$-adic norm 
$\|\cdot\|_w\,:\,w^*L_{\l}\ra \R$ and we see that the family 
$\|\cdot\|_v:=(\|\cdot\|_w)_{w\in W(F_v)}$ is a $v$-adic metric
on $L_{\l}$.
The family $(\|\cdot\|_v)_{v\in {\rm Val}(F)}$ 
will then be an adelic metric
on $L_{\l}$ (cf. \cite{Peyre} for $\l=2\rho_P$ and 
\cite{ST} for the definitions of
``$v$-adic metric''  and  ``adelic metric'').
The metrized line
bundle $\left(L_{\l}, (\|\cdot\|_v)_{v}\right)$ 
will be denoted by ${\cal L}_{\l}$.

\subsection*\no
{\bf 2.3}\hskip 0,5cm
Define a map
$$
H_P=H_{P,{\bf K}_G}\,:\, G({\bf A})\ra \Hom_{\C}(X^*(P)_{\C},\C)
$$
by $\langle\l,H_P(g)\rangle=\log (\prod_v|\l(p_v)|_v)$ for 
$g=pk$ with $ p=(p_v)_v\in P({\bf A}), k\in {\bf K}_G$
and $ \l\in X^*(P)$.
For $w=\pi_W(\g)\in W(F)$ and $ \g\in G(F)$ a simple computation (\cite{FMT})
shows that 
$$
H_{{\cal L}_{\l}}(w)=e^{-\langle \l,H_P(\g)\rangle}.
$$
The height zeta function 
$$
Z_W({\cal L}_{\l}, s)=\sum_{w\in W(F)}H_{{\cal L}_{\l}}(w)^{-s}
$$
is therefore an Eisenstein series 
$$
E^G_P(s\l-\rho_P,1_G)=
\sum_{\g\in P(F)\backslash G(F)}e^{\langle s\l,H_P(\g)\rangle}.
$$
To describe the domain of absolute convergence 
of this series we let $\Delta_0$ 
be the basis of positive roots of the root system $\Phi(S_0,G)$ which is
determined by $P_0$. For any $\a \in \Delta_0$ denote by $\check{\a}$ the
corresponding coroot. 
For $\l\in X^*(P)=X^*(S_0)$ we define $\langle\l,\a\rangle$ by
$(\l\circ\check{\a})(t)=t^{\langle\l,\a\rangle}$ and extend
this linearly in $\l$ to $X^*(P_0)_{\C}$.  Restriction of 
characters defines an inclusion $X^*(P)\ra X^*(P_0)$. Let 
$$
\Delta_0^P=\{\,\a\in \Delta_0\,\,|\,\, \langle\,\cdot\,,\a\rangle\,\, 
{\rm vanishes}\,\,{\rm on}\,\, X^*(P)\,\}, \,\,\,\,\, 
\D_P=\Delta_0 - \Delta_0^P.
$$
Put 
$$
X^*(P)^{+}=\{\l\in X^*(P)_{\R}\,\,|\,\,
\langle\l,\a \rangle > 0 \,\, {\rm for}\hskip 0,3cm
{\rm all}\,\,\, \a\in \Delta_P\,\}.
$$
By \cite{G}, Th\'eor\`eme 3, the Eisenstein series
$$
E^G_P(\l,g)=\sum_{\g \in P(F)\backslash G(F)}e^{
\langle\l+\rho_P,H_P(\g g)\rangle}
$$
converges absolutely for $\Re (\l)-\rho_P$ in 
$X^*(P)^+ $ and it can be meromorphically
continued to $X^*(P)_{\C}$ (cf. \cite{MW}, IV, 1.8).
The closure of the 
image of $X^*(P)^{+}$ in ${\rm Pic}(W)_{\R}$ is the cone 
$\L_{\rm eff}(W)$ generated by the effective divisors on $W$
(\cite{J}, II, 2.6).

\section{Toric varieties}
\label{3}

\subsection*\no
{\bf 3.1}\hskip 0,5cm
Let $T$ be a split algebraic torus of dimension $d$ over $F$. 
We put $M=X^*(T)$ and $N=\Hom(M,\Z)$. Let $\Sigma$ be a complete regular fan in
$N_{\R}$ such that the corresponding smooth toric variety $X=X_{\Sigma}$ is 
projective (cf. \cite{BaTschi1, Oda}). 
The variety $X$ is covered by affine open sets
$$
U_{\s}={\rm Spec}(F[M\cap \check{\s}]),
$$
where $\s$ runs through $\Sigma$ and $\check{\s}$ is the dual cone
$$
\check{\s}=\{\,m\in M_{\R}\,|\, n(m)\ge 0\,\, \forall\, n\in \s\,\}. 
$$
Denote by $PL(\Sigma)$ the group of $\Sigma$-piecewise linear integral 
functions on $N_{\R}$. By definition, a function $\p\,:\, N_{\R}\ra \R$
belongs to $PL(\Sigma)$ if and only if $\p(N)\in \Z$ and the restriction of
$\p$ to every $\s\in \Sigma$ is the restriction to $\s$ of a linear function 
on $N_{\R}$. 
For $\p\in PL(\Sigma)$ and every $d$-dimensional cone $\s\in \Sigma$ there exists
a unique $m_{\p,\s}\in M$ such that  for all $n\in \s$
we have
$$
\p(n)=n(m_{\p,\s}).
$$
Fixing for any $\s \in \Sigma$ a $d$-dimensional cone $\s'$ containing $\s$ we
put
$$
m_{\p,\s}=m_{\p,\s'}.
$$
To any $\p\in PL(\Sigma)$ we associate an invertible sheaf $L_{\p}$ on $X$
as the subsheaf of rational functions on $X$ generated over
$U_{\s}$ by $\frac{1}{m_{\p,\s}}$, considered as a rational function
on $X$ ($L_{\p}$ does not depend on the choice made above).
The assignment $\p\mapsto L_{\p}$ gives an exact sequence 
$$
0\ra M\ra PL(\Sigma)\ra \Pic (X)\ra 0
$$
(cf. \cite{Oda}, Corollary 2.5).

Denote by $\theta \,:\, X\times T\ra X$ the action of 
$T$ on $X$ and by $p_1\,:\, X\times T\ra X$ 
the projection  onto the first factor.
The induced $T$-action on the sheaf of rational functions restricts to
any subsheaf $L_{\p}$, i.e., there is a canonical $T$-linearization 
$$
\theta_{\p}\,:\, \theta^*L_{\p}\ra p_1^*L_{\p}
$$
(cf. \cite{MFK}, Ch. 1, \S\, 3, for the notion of a $T$-linearization).
In section four we will always consider $L_{\p}$ not merely as a line bundle
on $X$ but as a $T$-linearized line bundle with this $T$-linearization.
In this sense $PL(\Sigma)$ is isomorphic to the group 
$\Pic^T(X)$ of isomorphism 
classes of $T$-linearized line bundles on $X$. 

Let $\Sigma_1\subset N$ be the set of primitive integral generators of the 
one-dimensional cones in $\Sigma$ and put
$$
PL(\Sigma)^+:=\{\p\in PL(\Sigma)_{\R}\,\, |\,\, \p(e)> 0\,\, {\rm for}
\hskip 0,3cm {\rm all}\,\,\, e\in \Sigma_1\}.
$$
It is well-known (cf. \cite{Reid}, \cite{BaTschi1} Prop. 1.2.11),
that the cone of effective divisors 
$\L_{\rm eff}(X)\subset \Pic(X)_{\R}$ is the closure of the 
image of $PL(\Sigma)^+$ under the projection 
$PL(\Sigma)_{\R}\ra \Pic(X)_{\R}$.
Further, the anti-canonical line bundle on 
$X$ is isomorphic to $L_{\p_{\Sigma}}$,
where 
$\p_{\Sigma}(e)=1$ for all $e\in \Sigma_1$ 
(cf. \cite{BaTschi1}, Prop. 1.2.12).

\subsection*
\no
{\bf 3.2}\hskip 0,5cm
We shall introduce an adelic metric on the line bundle $L_{\p}$ as follows.
For 
$\s\in \Sigma$ and $v\in \Val (F)$ define
$$
{\bf K}_{\s,v}:=\{\,x\in U_{\s}(F_v)\,|\, |m(x)|_v\le 1\,\,\,
 \forall\,\, m\in \check{\s}\cap M\,\}.
$$These subsets cover $X(F_v)$ and we put for $x\in {\bf K}_{\s,v}$ and any
local section $s$ of $L_{\p}$ at $x$
$$
\|x^*s\|_x:=|s(x)m_{\p,\s}(x)|_v.
$$
The family $\|\cdot\|_v=(\|\cdot\|_x)_{x\in X(F_v)}$ is then a $v$-adic metric
on $L_{\p}$ and ${\cal L}_{\p}=\left(L_{\p},(\|\cdot\|_v)_v\right)$ 
is a metrization 
of $L_{\p}$.
Let ${\bf K}_{T,v}\subset T(F_v)$ be the maximal compact subgroup. 
Assigning to $x\in T(F_v)$ the map 
$$
M\ra \Z \,\,({\rm resp.}\,\, \R\,\, {\rm if}\,\, v|\infty),
$$
$$
m\mapsto -\frac{\log(|m(x)|_v)}{\log(q_v)},
$$
(where $q_v$ is the order of the residue field 
of $F_v$ for non-archimedean valuations
and $\log(q_v)=1$ for archimedean valuations)
we get an isomorphism 
$T(F_v)/{\bf K}_{T,v}\ra N$ (resp. $N_{\R}$ if $v|\infty$).
We will denote by $\overline{x}$ the image of 
$x\in T(F_v)$ in $N$ (resp. $N_{\R}$). For $\p\in PL(\Sigma)_{\C}$ define
$$
H_{\Sigma,v}(\,\cdot\,, \p)\,:\, T(F_v)\ra \C,
$$
$$
H_{\Sigma,v}(x,\p):= e^{-\p(\overline{x})\log (q_v)}.
$$
The corresponding global function  
$H_{\Sigma}(\,\cdot\,, \p)\, :\, T({\bf A})\ra \C$,
$$
H_{\Sigma}(x, \p):=\prod_v H_{\Sigma,v}(x_v,\p), 
$$
is well defined since for almost all $v$ the local component
$x_v$  belongs to
${\bf K}_{T,v}$.
The functions $H_{\Sigma,v}(\,\cdot\,, \p)$, $\p\in PL(\Sigma)$,
are related to the 
$v$-adic metric on $L_{\p}$ by the identity

\begin{equation}
\label{3.2.1}
H_{\Sigma,v}(x,\p)=\|x^*s_{\p}\|_x, \hskip 0,5cm (x\in T(F_v)),
\end{equation}
where $s_{\p}\in H^0(T,L_{\p})$ is the constant function $1$.
In particular, for every $x\in T(F)$ we have
$$
H_{{\cal L}_{\p}}(x)=H_{\Sigma}(x,\p)^{-1}.
$$

\subsection*\no
{\bf 3.3}\hskip 0,5cm
Let ${\bf K}_T=\prod_v{\bf K}_{T,v}\subset T({\bf A})$, and denote by
$$
{\cal A}_T=(T({\bf A})/T(F){\bf K}_T)^*
$$
the group of unitary characters of $T({\bf A})$ which are trivial on 
the closed subgroup $T(F){\bf K}_T$.
For $m\in M$ we 
obtain characters  $\c^m$ defined by
$$
\c^m(x):= e^{i\log(|m(x)|_{\bf A})}.
$$
This gives an embedding $M_{\R}\ra {\cal A}_T$. 
For any archimedean place $v$ and 
$\chi\in {\cal A}_T$ there is an $m_v=m_v(\chi)\in M_{\R}$ such that
$\chi_v(x_v)=e^{-i\overline{x}_v(m_v)}$ for all $x_v\in T(F_v)$. We
get a homomorphism 
\begin{equation}
\label{3.3.1}
{\cal A}_T\ra M_{\R,\infty}=\oplus_{v|\infty}M_{\R},
\end{equation}
$$
\chi\mapsto m_{\infty}(\chi)=(m_v(\chi))_{v|\infty}.
$$
Define $T({\bf A})^1$ to be the kernel of all maps $T({\bf A})\ra \R_{>0},
x\mapsto |m(x)|_{\bf A}$, for $m\in M$, and put
$$
{\cal U}_T=(T({\bf A})^1/T(F){\bf K}_T)^*.
$$
The choice of a projection ${\bf G}_m({\bf A})\ra {\bf G}_m({\bf A})^1$
induces by means of an isomorphism 
$T\stackrel{\sim}\longrightarrow
{\bf G}_{m,F}^d$  a splitting of the exact sequence
$$
1\ra T({\bf A})^1\ra T({\bf A})\ra T({\bf A})/T({\bf A})^1\ra 1.
$$
This gives decompositions
\begin{equation}
\label{3.3.2}
{\cal A}_T=M_{\R}\oplus {\cal U}_T
\end{equation}
and
$$
M_{\R,\infty}=M_{\R}\oplus M^1_{\R,\infty},
$$
where $M^1_{\R,\infty}$ is the minimal 
$\R$-subspace of $M_{\R,\infty}$ containing
the image of ${\cal U}_T$ under the map (\ref{3.3.1}). 
>From now on we fix such a (non-canonical) splitting. By
Dirichlet's unit theorem, the image of ${\cal U}_T\ra M^1_{\R,\infty}$ is 
a lattice of maximal rank. Its kernel is isomorphic to the character
group of ${\rm Cl}_F^d$, where ${\rm Cl}_F$ is the ideal class group of $F$.

For finite $v$ we let $dx_v$ be the Haar measure on $T(F_v)$ giving 
${\bf K}_{T,v}$ the volume one. 
For archimedean $v$ we take on $T(F_v)/{\bf K}_{T,v}$ the pull-back of
the Lebesgue measure on $N_{\R}$ (normalized by the lattice $N$) and on 
${\bf K}_{T,v}$ the Haar measure with total mass one. The product measure
gives an invariant measure  $dx_v$ on $T(F_v)$. On $T({\bf A})$ we get
a Haar measure $dx=\prod_vdx_v$.

\subsection*\no 
{\bf 3.4}\hskip 0,5cm 
We will denote by $S^1$ the unit circle. 
For a character $\chi\,:\, T(F_v)\ra S^1$ we define 
the Fourier transform of $H_{\Sigma,v}(\,\cdot\,, \p)$ by
$$
\hat{H}_{\Sigma,v}(\chi,\p)=\int_{T(F_v)}H_{\Sigma,v}(x_v,\p)\chi(x_v)dx_v.
$$
If $\chi$ is not 
trivial on ${\bf K}_{T,v}$ then $\hat{H}_{\Sigma,v}(\chi,\p)=0$
(assuming the convergence of the integral). 
We will show that these integrals do exists if 
$\Re (\p) $ is in $PL(\Sigma)^+$.

Let $v$ be an archimedean place of $F$. Any $d$-dimensional cone
$\s\in \Sigma$ is simplicial (since $\Sigma$ is regular) and it is generated by the set
$\s\cap \Sigma_1$. Let $\chi$ be unramified, i.e., 
$\chi(x)=e^{-i\overline{x}(m)}$ with some $m\in M_{\R}$. 
Then we get
\begin{equation}
\label{3.4.1}
\hat{H}_{\Sigma,v}(\chi,\p)=\sum_{\dim \s =d}\int_{\s}e^{-(\p(n)+in(m))}dn=
\sum_{\dim \s =d}\,\,\,\prod_{e\in \s\cap \Sigma_1}\frac{1}{\p(e)+ie(m)}.
\end{equation}
To give the result for finite places we define rational functions
$R_{\s}$ in variables $u_e, e\in \Sigma_1$, for any $\s\in \Sigma$ by
$$
R_{\s}((u_e)_e)=\prod_{e\in \s\cap \Sigma_1}\frac{u_e}{1-u_e},
$$
and put
$$
R_{\Sigma}((u_e)_e)=\sum_{\s\in \Sigma}R_{\s}((u_e)_e),
$$
$$
Q_{\Sigma}((u_e)_e)=R_{\Sigma}((u_e)_e)\prod_{e\in \Sigma_1}(1-u_e).
$$
Although elementary, it is a very important observation that the
polynomial $Q_{\Sigma}-1$ is a sum of monomials of degree not less than two
(cf. \cite{BaTschi1}, Prop. 2.2.3).

Let $\chi$ be an unramified unitary character of $T(F_v)$ and let 
$\Re(\p)$ be in $PL(\Sigma)^+$. Then we can calculate 
\begin{equation}
\label{3.4.2}
\hat{H}_{\Sigma,v}(\chi,\p)=\int_{T(F_v)}H_{\Sigma,v}(x_v,\p)\chi(x_v)dx_v=
\sum_{n\in N}e^{-\p(n)\log(q_v)}\chi(n)
\end{equation}
$$
= \sum_{\s\in \Sigma}\sum_{n\in {\s}^{o}\cap N}q_v^{-\p(n)}\chi(n)
=\sum_{\s\in \Sigma}R_{\s}\left((\chi(e)q_v^{-\p(e)})_e\right)
$$
$$
=Q_{\Sigma}\left((\chi(e)q_v^{-\p(e)})_e\right)
\prod_{e\in \Sigma_1}\left(1-\chi(e)q_v^{-\p(e)}\right)^{-1}.
$$
(Here we denoted by $\s^{o}$ the relative interior of the cone $\s$.)

Any  $e\in \Sigma_1$ induces a homomorphism $F[M]\ra F[\Z]$ and hence
a morphism of tori ${\bf G}_m\ra T$. For any character 
$\chi\in {\cal A}_T$ we denote by $\chi_e$ 
the Hecke character
$$
{\bf G}_m({\bf A})\longrightarrow
T({\bf A})\stackrel{\chi}\longrightarrow S^1
$$
thus obtained. 
The finite part of the Hecke 
$L$-function with character $\chi_e$ is by definition 
$$
L_f(\c_e,s)=\prod_{v\nmid \infty}(1-\chi_e(\pi_v)q_v^s)^{-1}
$$
and this product converges for $\Re (s)>1$ (here $\pi_v$ denotes
a local uniformizing element).

\no
By (\ref{3.4.1}) and (\ref{3.4.2}) we know that the global Fourier
transform 
$$
\hat{H}_{\Sigma}(\chi,\p)=\int_{T({\bf A})}H_{\Sigma}(x,\p)\chi(x)dx
$$
exists (i.e., the integral on the right converges absolutely) if
$\Re(\p)$ is contained in 
$\p_{\Sigma}+PL(\Sigma)^+$, because 
$$
\prod_{v\nmid\infty }
Q_{\Sigma}((\chi_v(e)q_v^{-\p(e)})_e)
$$ 
is an absolutely convergent
Euler product for $\Re (\p(e))>1/2 $ (for all $e\in \Sigma_1$) and hence
is bounded for 
$\Re (\p)$ in any compact subset in
$\frac{1}{2}\p_{\Sigma}+PL(\Sigma)^+$ (by some constant depending only
on this subset).

\bigskip
\no
{\bf Proposition 3.4}\,\,\,
{\it
The series 
$$
\sum_{x\in T(F)}H_{\Sigma}(xt,\p)
$$
converges absolutely and uniformly for $(\Re(\p),t)$ contained in any
compact subset of 
$(\p_{\Sigma}+PL(\Sigma)^+)\times T({\bf A}).$  
}

\bigskip
{\em Proof.}
Let ${\bf K}$ be a compact subset of $\p_{\Sigma}+PL(\Sigma)^+$ and let
$C_v\subset T(F_v)$ (for every $v\in \Val (F)$) be a compact subset,
equal to ${\bf K}_{T,v}$ for almost all $v$. 
Since any $\p\in PL(\Sigma)_{\C}$
is a continuous piecewise linear function 
(with respect to a finite subdivision of $N_{\R}$ into simplicial cones)
there exists a constant $c_v\ge 1$ (depending on ${\bf K}$ and $C_v$) 
such that for all $\p$ with $\Re(\p)\in {\bf K}$, 
$x_v\in T(F_v)$ and $t_v\in C_v$ we have
$$
\frac{1}{c_v}\le \left| 
\frac{H_{\Sigma,v}(x_vt_v,\p)}{H_{\Sigma,v}(x_v,\p)}\right|= 
\frac{H_{\Sigma,v}(x_vt_v,\Re(\p))}{H_{\Sigma,v}(x_v,\Re(\p))} \le c_v.
$$
If $C_v={\bf K}_{T,v}$ we may assume $c_v=1$. Put  $c=\prod_vc_v$.
For all $\p$ with $\Re(\p)\in {\bf K}$ and $t\in C:=\prod_v C_v$ we 
can estimate 
$$
\left| \sum_{x\in T(F)}H_{\Sigma}(xt,\p)\right|\le c\sum_{x\in T(F)}
H_{\Sigma}(x,\Re(\p)).
$$
Let $S$ be a finite set of places containing ${\rm Val}_{\infty}(F)$ and 
let $U_v\subset T(F_v)$ be a relatively compact open subset of
$T(F_v)$ for each $v\in S$, such that for all
$x_1\neq x_2\in T(F)$
$$
x_1U\cap x_2U=\emptyset,
$$
where $U=\prod_{v\in S}U_v\prod_{v\notin S}{\bf K}_{T,v}.
$
By the preceding argument, there exists a $c'>0$ such that for all
$\p\in {\bf K}$, $x\in T(F)$  and $u\in U$ 
$$
H_{\Sigma}(x,\p)\le c'H_{\Sigma}(xu,\p).
$$
Therefore,
$$
\sum_{x\in T(F)}H_{\Sigma}(x,\p)\le \frac{c'}{{\rm vol}(U)}\sum_{x\in T(F)}
\int_{U}H_{\Sigma}(xu,\p)du
$$
$$
\le \frac{c'}{{\rm vol}(U)}\int_{T({\bf A})}H_{\Sigma}(x,\p)dx <\infty
$$
by the discussion above. From the explicit expression for the integral
(cf. (\ref{3.5.1})) we derive the uniform convergence in $\p$ on ${\bf K}$. 
\hfill $\Box$

\subsection*
\no
{\bf 3.5}\hskip 0,5cm
The aim is to apply Poisson's summation formula to the height zeta
function. It remains to show that 
$\hat{H}(\,\cdot\,,\p)$ is absolutely integrable over ${\cal A}_T$.
For $\chi\in {\cal A}_T$ and $\Re (\p)$ contained in 
$\frac{1}{2}\p_{\Sigma}+PL(\Sigma)^+$ we put
$$
\zeta_{\Sigma}(\chi,\p):= \prod_{v|\infty}\hat{H}_{\Sigma,v}(\chi_v,\p)
\prod_{v\nmid\infty} Q_{\Sigma}((\chi_v(e)q_v^{-\p(e)})_e).
$$

\bigskip
\no
{\bf Lemma 3.5}\hskip 0,5cm
{\it
Let ${\bf K}$ be a compact subset of $PL(\Sigma)_{\C}$ such  that for all
$\p \in {\bf K}$ and $e\in \Sigma_1$
$$
\Re (\p(e))>\frac{1}{2}.
$$
Then there is a constant $c=c({\bf K})$ such that for all $\p\in {\bf K},
\chi\in {\cal A}_T$ and $m\in M_{\R}$ we have
$$
|\zeta_{\Sigma}(\chi,\p+im)|\le
 c\prod_{v|\infty}\left\{\sum_{\dim \s =d}\,\,\,
\prod_{e\in \s\cap \Sigma_1}\frac{1}{(1+
|e(m+m_v(\chi))|)^{1+1/d}}\right\}.
$$
}

{\em Proof.}
For ${\bf K}$ as above there exists a $c'>0$ such that for all $\chi\in
{\cal A}_T$ and $m\in M_{\R}$ one has
$$
\left|\prod_{v\nmid \infty}
Q_{\Sigma}((\chi_v(e)q_v^{-(\p(e)+ie(m))})_e)\right|\le c'
$$
for all $\p\in {\bf K}$  (see the argument before Proposition 3.4.
By \cite{BaTschi1}, Prop. 2.3.2, for all $v \mid {\infty}$
there is a constant $c_v$ such that for all 
$\p\in {\bf K}$, $\c\in {\cal A}_T$
and $m\in M_{\R}$ 
$$
|\hat{H}_{\Sigma,v}(\chi_v,\p+im)|\le c_v\sum_{\dim \s =d}
\,\,\,\prod_{e\in \s\cap \Sigma_1}
\frac{1}{(1+|e(m+m_v(\chi))|)^{1+1/d}}.
$$
Putting $c=c'\prod_{v|\infty}c_v$ we get the result.
\hfill $\Box$

\bigskip
\no
For $\Re (\p)$ contained in $\p_{\Sigma}+PL(\Sigma)^+$ 
we can write

\begin{equation}
\label{3.5.1}
\hat{H}_{\Sigma}(\chi,\p)=
\zeta_{\Sigma}(\chi,\p)\prod_{e\in \Sigma_1}L_f(\chi_e,\p(e)).
\end{equation}

\no
By the preceding lemma, we see that $\hat{H}_{\Sigma}(\,\cdot\,, \p)$ is 
absolutely integrable over ${\cal A}_T$. 
For $t\in T({\bf A})$ we have
$$
\int_{T({\bf A}_F)}H_{\Sigma}(xt,\p)
\chi(x)dx=\chi^{-1}(t)\hat{H}_{\Sigma}(\chi,\p).
$$
Hence we can apply Poisson's summation formula (together with (\ref{3.3.2}))
and obtain 
\begin{equation}
\label{3.5.2}
\sum_{x\in T(F)}H_{\Sigma}(xt,\p)= \mu_T\int_{M_{\R}}
\left\{\sum_{\c\in {\cal U}_T}
\hat{H}_{\Sigma}(\c,\p+im)(\c\c^m(t))^{-1}\right\}dm,
\end{equation}
where  the Lebesgue measure $dm$ on $M_{\R}$
is normalized by $M$ and 
$$
\mu_T=\frac{1}{(2\pi \kappa)^d},\hskip 0,5cm
\kappa=\frac{{\rm cl}_F\cdot {\rm R}_F}{{\rm w}_F}
$$ 
with ${\rm cl}_F$ the class number, ${\rm R}_F$ the regulator and 
${\rm w}_F$ the number of roots of unity in $F$.
Note that $\hat{H}_{\Sigma}(\c\cdot\c^m,\p)=\hat{H}_{\Sigma}(\c,\p+im).$

\subsection*\no
{\bf 3.6}\hskip 0,5cm
In section 5 we need uniform estimates for $L$-functions in a 
neighborhood of the line $\Re(s)=1$. For any unramified character
$\chi\,:\, {\bf G}_m({\bf A})/{\bf G}_m(F)\ra S^1$ and any
archimedean place $v$ there exists a $\tau_v\in \R$ such that
$\chi_v(x_v)=|x_v|_v^{i\tau_v}$ for all $x_v\in {\bf G}_m(F_v)$.
We put 
$$
\chi_{\infty}=(\tau_v)_{v|\infty}\in \R^{{\rm Val}_{\infty}(F)}
\hskip 0,5cm {\rm and}\hskip 0,5cm 
\|\chi_{\infty}\|=\max_{v|\infty}|\tau_v|.
$$
We will use
the following theorem of Rademacher (\cite{Rademacher}, Theorems 4,5),
which rests on the Phragm\'en-Lindel\"of principle. 

\bigskip
\no
{\bf Theorem 3.6}\hskip 0,5cm {\it
For any $\epsilon>0$ there exists a $\d>0$ and a constant $c(\epsilon)>0$
such that
for all $s$ with $\Re(s)> 1-\d$ and all unramified Hecke characters
which are non-trivial on ${\bf G}_m({\bf A})^1$ one has
\begin{equation}
\label{3.6.1}
|L_f(\chi,s)|\le c(\epsilon)(1+|\Im (s)| +\|\chi_{\infty}\|)^{\epsilon}.
\end{equation}
For the trivial character $\chi=1$ one has
\begin{equation}
\label{3.6.2}
|L_f(1,s)|\le c(\epsilon)\left|\frac{1+s}{1-s}\right|(1+|\Im(s)|)^{\epsilon}.
\end{equation}
}

\bigskip
\bigskip

\section{Twisted products}
\label{4}

\subsection*\no
{\bf 4.1}\hskip 0,5cm
Let $G,P,W=P\backslash G$ etc. be as in section 2 and 
$T,\Sigma,X=X_{\Sigma}$ etc.
be as in section 3. 
Let $\eta \,:\, P\ra T$ be a homomorphism. Then $P$ acts from the right on 
$X\times G$ by 
$$
(x,g)\cdot p := (x\eta (p),p^{-1}g).
$$
Since $\pi_W\,:\, G\ra W$ is locally trivial, the quotient
$$
Y=X\times^{P} G:= (X\times G)/P
$$
exists as a variety over $F$. Moreover, the projection 
$X\times G\ra G$ induces a morphism $\pi \,: \, Y\ra W$ 
and $Y$ becomes a locally trivial fiber bundle over $W$ with fiber $X$ 
(compare \cite{J}, I.5.16). Hence, by the properties 
of $X$ (non-singular, projective), we 
see that $Y$ is a non-singular projective variety over $F$ (``projectivity''
requires a short argument, cf. \cite{Strauch}). 
The quotient morphism $X\times G\ra Y$ will be denoted by $\pi_Y$.
Let $\p\in PL(\Sigma)$ and let $L_{\p}$ be the invertible sheaf on $X$ defined
in section 3.1. Denote by ${\bf L}_{\p}$ the corresponding 
${\bf G}_a$-bundle over $X$, i.e., ${\bf L}_{\p}={\bf V}(L^{\vee}_{\p})=
{\bf V}(L_{-\p})$ (with the notation of \cite{Ha}, II, Exercise 5.18).
The canonical $T$-linearization
$$
\theta_{-\p}\,:\, \theta^*L_{-\p}\ra p_1^*L_{-\p}
$$
induces an action ${\bf L}_{\p}\times T\ra {\bf L}_{\p}$ of $T$ on
${\bf L}_{\p}$ which is compatible with the action of $T$ on $X$. 
The twisted product ${\bf L}_{\p}\times^PG$ will then be a ${\bf G}_a$-bundle
over $Y$ and we define $L^Y_{\p}$ to be the sheaf of local sections of
${\bf L}_{\p}\times^PG $ over $Y$.  Note that $L^Y_{\p}$ (and even
its isomorphism class in ${\rm Pic}(Y) $) 
depends on the fixed $T$-linearization
$\theta_{-\p}$. In fact, for $\p\in PL(\Sigma)$ and $m\in M$ we have
$$
L^Y_{\p +m}\simeq L^Y_{\p}\otimes \pi^* L_{m\circ \eta}.
$$
Embedding $M$ in $PL(\Sigma)\oplus X^*(P)$ by $m\mapsto (m,-m\circ\eta)$ 
we see that $M$ is contained in the kernel of the homomorphism 
$$
\begin{array}{ccl}
\psi \,:\, PL(\Sigma)\times X^*(P) &\ra &{\rm Pic}(Y),\\
(\p,\l) &\mapsto & {\rm isomorphism}\,\, {\rm class}\,\, {\rm of}\,\, 
L^Y_{\p}\otimes \pi^*L_{\l}.
\end{array}
$$

\subsection*\no
{\bf 4.2}\hskip 0,5cm
In the following proposition we collect 
all relevant facts about the geometry
of twisted products which we will need in the sequel.

\bigskip
\no
{\bf Proposition 4.2}\hskip 0,5cm
{\it
a) The sequence 
$$
0\ra M\ra PL(\Sigma)\oplus X^*(P)\ra {\rm Pic}(Y)\ra 0
$$
is exact.

b) The cone of effective divisors $\L_{\rm eff}(Y)\subset {\rm Pic}(Y)_{\R}$
is the image of the closure of
$$
PL(\Sigma)^{+}\times X^*(P)^{+}\subset PL(\Sigma)_{\R}\oplus X^*(P)_{\R}.
$$

c) The anti-canonical line bundle $\KX_Y$ is isomorphic to
$L^Y_{\p_{\Sigma}}\otimes \pi^*L_{2\rho_P}$.
}
\bigskip

{\em Proof.}
a) By \cite{Sa}, Proposition 6.10, there is an exact sequence
$$
F[X\times G]^*/F^*\ra X^*(P)\ra \Pic(Y)\ra\Pic(X\times G). 
$$
Denote by $\pi_X\,:\, X\times G\ra X$ the canonical projection. Let $L$ be 
an invertible sheaf on $Y$. Then 
$$
\pi^*_YL\simeq \pi^*_XL_{\p}
$$
(for some $\p\in PL(\Sigma)$) because $\Pic(X\times G)=
{\rm Pic}(X)\oplus \Pic(G)=\Pic(X)$ 
(cf. \cite{Sa}, Lemme 6.6 (i) and Lemme 6.9 (iv)).
Note that $\pi^*_YL^Y_{\p}\simeq \pi^*_XL_{\p}$, so that 
$$
\pi^*_Y(L\otimes L^Y_{-\p})
$$
is trivial. Hence there exists a character $\l $ of $P$ such that
$L\otimes L^Y_{-\p}$ is isomorphic to $\pi^*L_{\l} $
 (the map $X^*(P)\ra {\rm Pic}(Y)$ factorizes 
$X^*(P)\ra {\rm Pic}(W)\ra {\rm Pic}(Y)$). This shows surjectivity. 

Suppose now that for $\p \in PL(\Sigma)$ and $ \l \in X^*(P)$ 
the sheaf $L_{\p}^Y\otimes \pi^*L_{\l}$ is trivial
on $Y$. Then $\pi_Y^*(L^Y_{\p}\otimes \pi^*L_{\l})\simeq \pi^*_XL_{\p}$
is trivial on $X\times G$, therefore $L_{\p}\simeq {\cal O}_X, \p=m\in M$ 
and $L^Y_{\p}\otimes \pi^*L_{\l}=\pi^*L_{\l+m\circ\eta}$. By Rosenlicht's
theorem, 
$$
F[X\times G]^*/F^*=F[X]^*/F^*\oplus F[G]^*/F^*=X^*(G)=0, 
$$
therefore, 
the map $X^*(P)\ra {\rm Pic}(Y)$ is injective, 
hence $\l+m\circ \eta =0$ and $(\p,\l)$ is
in the image of $M\ra PL(\Sigma)\oplus X^*(P)$.

b) For $\p\in PL(\Sigma)$ denote by $\Box_{\p}$ the set of
 all $m\in M$ such that for all $n\in N_{\R}$
$$
\p(n)+n(m)\ge 0.
$$
By \cite{Oda}, Lemma 2.3, $\Box_{\p}$ is a basis for $H^0(X,L_{\p})$ 
(note the different sign conventions). It is easy to see that
$$
\pi_*L_{\p}^Y\simeq \oplus_{m\in \Box_{\p}}L_{-m\circ\eta}.
$$
Suppose $L^Y_{\p}\otimes \pi^*L_{\l}$ has a non-zero global section. Then
$$
\pi_*(L_{\p}^Y\otimes \pi^*L_{\l})\simeq
 \oplus_{m\in \Box_{\p}}L_{-m\circ \eta +\l}
$$
has a non-zero global section, hence (cf. section 2.3) there is a
$m'\in \Box_{\p}$ such that $-m'\circ \eta +\l$ is contained in the
closure of $ X^*(P)^{+}$. 
Putting $\l'=-m'\circ \eta +\l, \p'=\p+m'$ we have
$L^Y_{\p'}\otimes \pi^*L_{\l'}\simeq L^Y_{\p}\otimes \pi^*L_{\l}$ and
$(\p',\l')$ is contained in the closure of
$PL(\Sigma)^+\times X^*(P)^+$. 

On the other hand, if $(\p,\l)\in PL(\Sigma)\oplus X^*(P)$
is contained in the closure of $PL(\Sigma)^+\times X^*(P)^+$
then the trivial character corresponds to a global section of 
$L_{\p}$. Hence 
$$
\pi_*(L_{\p}^Y\otimes \pi^*L_{\l})
= L_{\l}\oplus \bigoplus_{m\in \Box_{\p}-\{0\}}L_{-m\circ \eta +\l}
$$
and $H^0(W,L_{\l})\neq \{0\}$, i.e., $L_{\p}^Y\otimes \pi^*L_{\l}$
has a non-zero global section. 

c) Note first that the exact sequence
$$
0\ra \pi^*\Omega_W\ra \Omega_Y\ra \Omega_{Y/W}\ra 0
$$
splits, 
and therefore $\omega_Y\simeq (\L^d\Omega_{Y/W})\otimes \pi^*\omega_W$.
Since $\omega_W\simeq L_{-2\rho}$ it remains to show 
that $(\L^d\Omega_{Y/W})^{\vee}
\simeq L^Y_{\p_{\Sigma}}$. Let ${\cal J}_{Y/W}$ be the 
ideal sheaf of the image of the diagonal morphism 
$$
\Delta_{Y/W}\,:\, Y\ra Y\times_WY.
$$
But $Y\times_WY$ is canonically isomorphic to $(X\times X)\times^PG$ and 
$\Delta_{Y/W}(Y)$ is just $\Delta_X(X)\times^PG$. 
Hence we see that $({\cal J}_{Y/W}/{\cal J}_{Y/W}^2)^{\vee}$ is the sheaf of
local sections of
$$
{\bf V}({\cal J}_X/{\cal J}_X^2)\times^PG,
$$
where ${\cal J}_X$ is the ideal sheaf of $\Delta_X(X)\subset X\times X$. 
Pulling back to $Y$ and taking the $d$-th exterior power we get
$$
{\bf V}(\L^d\Omega_{Y/W})\simeq {\bf V}(\L^d\Omega_X)\times^PG 
\simeq {\bf V}(\omega_X)\times^PG.
$$
The canonical $T$-linearization of $\omega_X$ (induced by the action of
$T$ on rational functions) corresponds to the 
$T$-linearization $\theta_{-\p_{\Sigma}}$ of 
$L_{-\p_{\Sigma}}\simeq \omega_X$, i.e.,
$$
{\bf L}_{\p_{\Sigma}}\times^PG\simeq {\bf V}(\omega_X)\times^PG
$$
and we get $L^Y_{\p_{\Sigma}}\simeq (\L^d\Omega_{Y/W})^{\vee}$.
\hfill $\Box$

\subsection*
\no
{\bf 4.3}\hskip 0,5cm
We are going to introduce an adelic metric on the sheaves $L^Y_{\p}$.
A section of $L^Y_{\p}$ over an open subset
$U\subset Y$ can be identified with  a $P$-equivariant morphism
$s\,:\, \pi^{-1}_Y(U)\ra {\bf L}_{\p}$ over 
$X$, i.e.,
$$
s(x\eta(p),p^{-1}g)=s(x,g)\cdot \eta(p).
$$
Let $v$ be a place of $F$ and let $y\in Y(F_v)$ be the image of
$(x,k)\in X(F_v)\times G(F_v)$ with $k\in {\bf K}_{G,v}$. Let
$s\,:\, \pi^{-1}_Y(U)\ra {\bf L}_{\p}$ be a local section of $L^Y_{\p}$
over $U\subset Y$ with $y\in U(F_v)$. Define
$$
\|\cdot\|_y\,:\, y^*L^Y_{\p}\ra \R
$$
by 
$$
\|y^*s\|_y=\|s\circ (x,k)\|_x.
$$
Then $\|\cdot \|_v=(\|\cdot \|)_{y\in Y(F_v)}$ is a $v$-adic metric on 
$L^Y_{\p}$ and 
${\cal L}^Y_{\p}=\left(L^Y_{\p},(\|\cdot \|_v)_v\right)$ is a metrization 
of $L^Y_{\p}$ (cf. \cite{Strauch}).
Let 
$$
Y^o=T\times^P G\hookrightarrow X\times^PG=Y
$$
be the twisted product of 
$T$ with $W$. Over $Y^o$ there is a canonical section of $L^Y_{\p}$, namely
$$
s^Y_{\p}\,:\, \pi^{-1}(Y^o)=T\times G\ra {\bf L}_{\p},
$$
$ s^Y_{\p}(x,g)=s_{\p}(x)$, where $s_{\p}\in H^0(T,L_{\p})$
corresponds to the constant function $1$. 
Let $y=\pi_W(x,g)\in Y(F_v)$ where $ g=pk$ with 
$ p\in P(F_v)$ and $ k\in {\bf K}_{G,v}$. 
Then 
$$
\|y^*s_{\p}^Y\|_y = \|s_{\p}^Y\circ (x\eta(p),k)\|_{x\eta(p)} 
$$
$$
= \|(x\eta(p))^*s_{\p}\|_{x\eta(p)}=e^{-\p(\overline{x\eta(p)})\log (q_v)}
$$
(by (\ref{3.2.1})).
Globally, for $y\in Y^o(F), y=\pi_Y(x,\g), x\in T(F), 
\g\in G(F),\g=p_{\g}k_{\g}$ and $ p_{\g}, k_{\g}$ as above, in
$P({\bf A}), {\bf K}_G$, respectively,  
we get 
\begin{equation}
\label{4.3.1}
H_{{\cal L}^Y_{\p}}(y)=
\prod_v\|y^*_vs^Y_{\p}\|_{y_v}^{-1}=H_{\Sigma}(x\eta(p_{\g}),-\p).
\end{equation}

\subsection*\no
{\bf 4.4}\hskip 0,5cm
Let $\xi \,:\, P({\bf A})/P(F)\ra S^1$ be an unramified character, i.e., 
$\xi$ is trivial on $P({\bf A})\cap {\bf K}_G$. 
Using the Iwasawa decomposition  we get a well defined function
$$
\phi_{\xi}\,:\, G({\bf A})\ra S^1,
$$
$$
\phi_{\xi}(g)=\xi(p),
$$
if $g=pk$ as above. We denote by 
$$
E^G_P(\l,\xi,g)=
\sum_{\g\in P(F)\backslash G(F)}\phi_{\xi}(\g g)
e^{\langle\l+\rho_P,H_P(\g g)\rangle}
$$
the corresponding Eisenstein series and we put
$E^G_P(\l,\xi)=E^G_P(\l,\xi,1_G)$. This series converges absolutely for
${\rm Re}(\l)$ contained in the cone $\rho_P+X^*(P)^+$
(cf. (2.3)). A character $\c\in {\cal A}_T$ induces a character
$\c_{\e}=\c\circ\e\,:\, P({\bf A})/P(F)\ra S^1$. We denote by
$\check{\eta}\,:\, X^*(T)_{\R}\ra X^*(P)_{\R}$ the map on characters
induced by $\eta$.

\bigskip
{\bf Proposition 4.4}
{\it
Let $L$ be a line bundle on $Y$ such that its class is contained in
the interior of the cone
$\L_{\rm eff}(Y)$. Let $(\p,\l)$ be in $PL(\Sigma)^+\times
X^*(P)^+$ with $\psi (\p,\l)=[L]$.
There is a metrization ${\cal L}$ of $L$ such that for all $s$ with
${\rm Re}(s)(\p,\l)\in (\p_{\Sigma},2\rho_P) + PL(\Sigma)^+\times X^*(P)^+$
the series 
$$
Z_{Y^o}({\cal L},s)=\sum_{y\in Y^{o}(F)}H_{{\cal L}}(y)^{-s}
$$
converges absolutely. Moreover, for these $s$
$$
Z_{Y^o}({\cal L},s)=\mu_T\int_{M_{\R}}
\left\{
\sum_{\c\in {\cal U}_T}\hat{H}_{\Sigma}(\c,s\p+im)
E^G_P(s\l-\rho_P - i\check{\eta}(m),\ce^{-1})
\right\}dm,
$$
where the sum and  integral on the right converge absolutely too. 
}
\bigskip

{\em Proof.}
Let $(\p',\l')\in PL(\Sigma)\oplus X^*(P)$ such that 
there is an isomorphism $L\simeq L^Y_{\p'}\oplus 
\pi^* L_{\l'}$. Denote by ${\cal L}$ the 
metrization of $L$ which is the pullback of 
${\cal L}_{\p'}^Y\otimes \pi^*{\cal L}_{\l'}$ via this isomorphism.
Let $m\in M_{\R}$ such that 
$$
(\p,\l)=(\p'+m,-\check{\eta}(m) +\l')
$$ 
is contained
in $PL(\Sigma)^+\times X^*(P)^+$. 
By (\ref{4.3.1}), we have
for any $y\in Y^o(F), y=\pi_Y(x,\g )$ with $ x\in T(F),$ $
\g \in G(F)$ and $ \g =p_{\g}k_{\g}$ 
$$
H_{{\cal L}}(y)=H_{{\cal L}_{\p'}^Y\otimes \pi^*{\cal L}_{\l'}}(y)=
e^{-\langle\l',H_P(\g )\rangle}H_{\Sigma}(x\eta (p_{\g}),-\p')
$$
$$
= e^{-\langle\l'-m\circ \eta,H_P(\g)\rangle}
H_{\Sigma}(x\eta(p_{\g}),-(\p'+m))= 
 e^{-\langle\l ,H_P(\g)\rangle}H_{\Sigma}(x\eta(p_{\g}),-\p).
$$
We consider $s=u+iv\in \C$ such that $u\cdot \p$ is contained in 
the shifted cone $\p_{\Sigma}+PL(\Sigma)^+$
and  $u\cdot \l$ is contained in the cone $2\rho_P+X^*(P)^+$. 
Then
$$
\sum_{x\in T(F)}H_{\Sigma}(x\eta(p_{\g}),u\p)
$$
converges by Proposition 3.4 and is equal to
$$
\mu_T\int_{M_{\R}}\{\sum_{\chi\in {\cal U}_T}
\hat{H}(\c\c^m,u\p)\c\c^m(\eta(p_{\g}))^{-1}\}dm
$$
(cf. (\ref{3.5.2})). Moreover, $\hat{H}_{\Sigma}(\,\cdot\, , u \p)$
is absolutely convergent on ${\cal A}_T$ and therefore
$$
\sum_{x\in T(F)}H_{\Sigma}(x\eta(p_{\g}), u\p)\le 
\mu_T\int_{M_{\R}}\left\{\sum_{\chi\in {\cal U}_T}
\left| \hat{H}_{\Sigma}(\chi\c^m,u\p)\right|\right\}dm 
$$
is bounded by some constant $c$ (which is 
independent of $\eta(p_{\g})$).
Thus we may calculate
$$
\sum_{y\in Y^o(F)}\left| H_{\cal L}(y)^{-s}\right|
=\sum_{\g\in P(F)\ba G(F)}
e^{\langle u\l,H_{P}(\g)\rangle}\sum_{x\in T(F)}H_{\Sigma}(x\eta(p_{\g}),
u\p)
$$
$$
\le c \sum_{\g\in P(F)\ba G(F)}e^{\langle u\l,H_{P}(\g)\rangle}.
$$
This shows the first assertion.
Since
$$
\mu_T\int_{M_\R}\left\{\sum_{\chi\in {\cal U}_T}
\left|\hat{H}_{\Sigma}(\chi\c^m, u\p)\right|\right\}dm
$$
converges, we can interchange the summation and integration  and get
$$
Z_{Y^o}({\cal L},s)
=\sum_{\g\in P(F)\backslash G(F)}
e^{\langle s\l,H_P(\g)\rangle}\mu_T\int_{M_{\R}}
\left\{\sum_{\c\in {\cal U}_T}\hat{H}_{\Sigma}(\c\c^m,s\p)(\c\c^m)^{-1}
(\eta(p_{\g}))\right\}dm
$$
$$
=\mu_T\int_{M_{\R}}\left\{\sum_{\chi\in {\cal U}_T}
\hat{H}_{\Sigma}(\chi,s\p+im)E^G_P(s\l-\rho_P-i\check{\eta}(m), 
\ce^{-1})\right\}dm.
$$
\hfill $\Box$

\bigskip
\bigskip

\section{Meromorphic continuation}
\label{5}

\subsection*\no
{\bf 5.1}\hskip 0,5cm
The proposition in section 4.4 gives an expression of the
height zeta function (for the open subset $Y^o\subset Y$ )
which we will use to determine the asymptotic behavior of the counting
function $N_{Y^o}({\cal L},H) $ (cf. sec. \ref{1}) by applying 
a Tauberian theorem.

The first thing to do is to show that $Z_{Y^o}({\cal L},s)$ can be
continued meromorphically to a half-space beyond the abscissa of 
convergence  and that there is no pole on this line
with non-zero imaginary part.
Then it remains to prove that this abscissa is at  $\Re(s)=a(L)$ and to
determine the order of the pole in $s=a(L)$. We will see that
this order is $b(L)$.

The method which we will explain now consists in an iterated application 
of Cauchy's residue theorem. The proofs will be given in section 7.

\subsection*
\no
{\bf 5.2}\hskip 0,5cm
Let $E$ be a finite dimensional vector space over $\R$ and $E_{\C}$
its complexification. Let $V\subset E$ be a subspace and let $l_1,...,l_m\in
E^{\vee} =\Hom_{\R}(E,\R)$
be linearly independent linear forms. Put 
$H_j={\rm Ker}(l_j)$ for $j=1,...,m$. 

Let $B\subset E$ be an 
open and convex neighborhood of ${\bf 0}$ such that 
for all $x\in B$ and $j=1,...,m$ we have $l_j(x)>-1$.
 Let $T_B=B+iE\subset E_{\C}$ be the tube domain over $B$ 
and denote by ${\cal M}(T_B)$ the set of meromorphic functions on $T_B$.
We consider meromorphic functions $f\in {\cal M}(T_B)$ with the following 
properties: The function 
$$
g(z)=f(z)\prod_{j=1}^m\frac{l_j(z)}{l_j(z)+1}
$$
is holomorphic in $T_B$ and there is a  sufficient function 
$c\,:\, V\ra \R_{\ge 0}$ such that for all compacts ${\bf K}\in T_B$, all
$z\in {\bf K}$ and all $v\in V$ we have the estimate
$$
|g(z+iv)|\le \kappa ({\bf K})c(v).
$$
(Cf. section 7.3 for a precise definition of a sufficient function. 
In particular, such a sufficient function is absolutely integrable over 
any subspace $U\subset V$.)
In this case we call $f$  {\it distinguished} with respect
to the data $(V;l_1,...,l_m)$. 

Let $C$ be a connected component of $B- \cup_{j=1}^mH_j$.
By the conditions on $g$ the integral
$$
\tilde{f}_C(z):=\frac{1}{(2\pi)^{\nu}}\int_V f(z+iv)dv
$$
($\nu =\dim V$ and $dv$ is a fixed Lebesgue measure on $V$) converges
for every $z\in T_C$ and $\tilde{f}_C$ is a holomorphic function on $T_C$.

\bigskip
\no
{\bf Theorem 5.2}\hskip 0,5cm
{\it
There is an open neighborhood $\tilde{B}$ containing $C$, and 
linear forms $\tilde{l}_1,...,\tilde{l}_{\tilde{m}}$ which vanish on $V$
such that
$$
z\mapsto \tilde{f}_{C}(z)\prod_{j=1}^{\tilde{m}}\tilde{l}_j(z)
$$
has a holomorphic continuation to $T_{\tilde{B}}$. Moreover, for all
$j\in \{1,...,\tilde{m}\}$ we have $\Ker (\tilde{l}_j)\cap C=\emptyset$.
}
\bigskip

\no
We shall give the proof of this theorem in sections 7.3 and 7.4.

\subsection*\no
{\bf 5.3}\hskip 0,5cm
Put $E^{(0)}=\cap_{j=1}^m\Ker (l_j)$ and $E_0=E/E^{(0)}$. Let
$\pi_0\,:\, E\ra E_0$ be the canonical projection and suppose 
$V\cap E^{(0)} = \{{\bf 0}\}$.  Let
$$
E^+_0=\{x\in E_0\,|\, l_j(x)\ge 0\,\, {\rm for}\,\,{\rm all} \,\, j=1,...,m\}
$$
and let $\psi_0\,:\, E_0\ra P:=E^+_0/\pi_0(V)$ be the canonical projection. We
want to assume that $\pi_0(V)\cap E^+_0=\{ {\bf 0}\}$, 
so that $\L:=\psi_0(E^+_0)$ is a strictly convex polyhedral cone. 
Let $dy$ be the  Lebesgue measure on $E_0^{\vee}$ normalized by the lattice
$\oplus_{j=1}^m\Z l_j$. Let $A\subset V$ be a lattice and let $dv$ be the
measure on $V$ normalized by $A$. On $V^{\vee}$ we have the Lebesgue 
measure $dy'$ normalized by $A^{\vee}$ and a section of the projection
$E^{\vee}_0\ra V^{\vee}$ gives a measure $dy''$ on $P^{\vee}$
with $dy=dy'dy''$.

Define the ${\cal X}$-function of the cone $\L$ by 
$$
{\cal X}_{\L}(x)=\int_{\L^{\vee}}e^{-y''(x)}dy''
$$
for all $x\in P_{\C}$ with $\Re(x)$ contained in the interior of
$\L$ (cf. section 7.1).

Let $B\subset E$ be as above and let $f\in {\cal M}(T_B)$ be a distinguished
function with respect to $(V;l_1,...,l_m)$. Put
$$
g(z)=f(z)\prod_{j=1}^m\frac{l_j(z)}{l_j(z)+1},
$$
$$
B^+=B\cap \{x\in E\,|\, l_j(x)\ge 0, \,{\rm for}\,\,
 {\rm all}\,\, j = 1,...,m\},
$$
$$
\tilde{f}_{B^+}(z)=\frac{1}{(2\pi )^{\nu}}\int_Vf(z+iv)dv
$$
(for $z\in T_B$).
The function $\tilde{f}_{B^+}\,:\, T_{B^+}\ra \C$ is holomorphic and has
a meromorphic continuation to a neighborhood of ${\bf 0}\in E_{\C}$.
In section 7.5 we will prove the following theorem.

\bigskip
\no
{\bf Theorem 5.3}\hskip 0,5cm {\it
For $x_0\in B^+$ we have
$$
\lim_{s\ra 0}s^{m-\nu}\tilde{f}_{B^+}(sx_0)=
g({\bf 0}){\cal X}_{\L}(\psi_0(x_0)).
$$
}

\bigskip

\subsection*\no
{\bf 5.4}\hskip 0,5cm
In this section we make some preparations in 
order to apply the general setting
of 5.2.
Let $L$ be a line bundle on $Y$ such that its class in ${\rm Pic}(Y)$ lies in 
the interior of $\L_{\rm eff}(Y)$. By the definition of $a(L)$ (cf. section 1
and Proposition 4.2),
$$
a(L)[L]-\psi(\varphi_{\Sigma},2\rho_P)\in \L(L)
$$
where $\L(L)$ is the minimal face of $\L_{\rm eff}(Y)$ containing
$a(L)[L]-\psi(\varphi_{\Sigma},2\rho_P)$. Define $\p_e\in PL(\Sigma) $ 
(for $e\in \Sigma_1$) by $\p_e(e')=\delta_{ee'}$, for all $e'\in \Sigma_1$
and put
$$
\Sigma'_1:= \{e\in \Sigma_1\,|\, \psi(\p_e,{\bf 0})\in \L(L)\}.
$$
Let $P'\subset G$ be the standard parabolic subgroup with
$$
\Delta_{P'}=\{\a\in \D_P\,\,|\,\, \psi({\bf 0},\w_{\a})\in \L(L)\},
$$ 
where $\langle\varpi_{\a},\b\rangle =\delta_{\a \b}$ for all
$\a ,\b\in \Delta_0$. 
Let 
$$
(\p_L,\l_L)\in (\sum_{e\in \Sigma'_1}\R_{>0}\p_e)\times 
(\sum_{\a\in \Delta_{P'}}\R_{>0}\w_{\a})
$$
such that $\psi(\p_L,\l_L)=a(L)[L]-\psi (\p_{\Sigma},2\rho_P)$. 
Then 
$$
\hat{L}:=\frac{1}{a(L)}(\p_{\Sigma}+\p_L,2\rho_P+\l_L)
$$
is mapped onto $[L]$ by $\psi$. 
Denote by 

\begin{equation}
\label{hL}
h_L(\p,\l):=\prod_{e\in \Sigma_1-\Sigma_1'}
\frac{\p(e)}{\p(e)+1}\prod_{\a\in \Delta_P - \Delta_{P'} }
\frac{\langle \l,\a\rangle}{
\langle \l,\a\rangle +1}
\end{equation}

\no 
and put
$$
\tp=\p+\p_{\Sigma}+ \p_L\hskip 0,5cm {\rm  and}
\hskip 0,5cm \tl=\l+\rho_P+\l_L.
$$

\bigskip
\no
>From now on we will denote by 
${\bf K}_G\subset G({\bf A})$ the maximal compact subgroup
defined in section 8.2.
\bigskip

\no
{\bf Lemma 5.4}\hskip 0,5cm
{\it
There exists a convex open neighborhood $B$ of ${\bf 0}$ in 
$PL(\Sigma)_{\R}\oplus X^*(P)_{\R}$ with the 
following property:
For any compact subset ${\bf K}\subset T_B$  
there is a constant $c=c({\bf K})>0$ such that for all
$(\p,\l)\in {\bf K}, \c\in {\cal U}_T$ and $m\in M_{\R}$
we have
$$
\left|\hat{H}_{\Sigma}(\chi,\tp +im)E^G_P(\tl-
i\check{\eta}(m),\c_{\eta}^{-1})h_L(\p+im,\l -i\check{\eta}(m))\right|
$$
$$
\le c\prod_{v|\infty}\left\{\sum_{\dim \s =d}\,\,\prod_{e\in \s\cap \Sigma_1}
\frac{1}{(1+|e(m+m_v(\chi))|)^{1+1/2d}}\right\}.
$$
}
\bigskip

{\em Proof.}
Write as in (\ref{3.5.1})
$$
\hat{H}_{\Sigma}(\c,\tp+im)=\zeta_{\Sigma}(\c,\tp+im)\prod_{e\in \Sigma_1}
L_f(\c_e,1+(\p+\p_L+im)(e)).
$$
For $\Re(\p)$ sufficiently small and $e\in \Sigma_1'$ we have
$$
\Re(\p(e))+\p_L(e)\ge \frac{1}{2}\p_L(e)>0.
$$
Hence
$$
\left|
\prod_{e\in \Sigma_1'}L_f(\c_e,1+(\p+\p_L+im)(e))
\right|
$$
is bounded for $\Re(\p)$ sufficiently small.
If $e\in \Sigma_1-\Sigma_1'$ then $\p_L(e)=0$. By the estimates of 
Rademacher (cf. Theorem 3.6), we have for $\c_e\neq 1$
$$
\left|
L_f(\c_e,1+(\p+im)(e))
\right|\le 
c_e(1+|m(e)|+\|(\c_e)_{\infty}\|)^{\epsilon}
$$
for  $\Re(\p(e))>-\delta$ and $\p $ in a compact set ($ \delta$
depends on $\epsilon $, $c_e$ depends on this compact subset).
\no
If $\c_e=1$ (abusing notations we will denote
from now on the trivial character by 1) then
$$
\frac{(\p+im)(e)}{(\p+im)(e)+1}
\left|
L_f(1,1+(\p+ im)(e))
\right|\le 
c_e(1+|m(e)|)^{\epsilon}
$$
Now we use Proposition 8.7 concerning estimates for Eisenstein
series. This proposition tells us that there is for given $\epsilon >0$ an
open neighborhood of ${\bf 0}$ in $X^*(P)_{\R}$ such that for $\Re(\l)$ 
contained in this neighborhood
$$
\left|
\prod_{\a\in \D_P}
\frac{\langle \l+\l_L-i\check{\eta}(m),\a\rangle}{
\langle \l+\l_L-i\check{\eta}(m),\a\rangle +1}
E^G_P(\tl -i\check{\eta}(m),\ce^{-1})
\right|
\le c_1(1+\|\Im(\l)+ 
\check{\eta}(m)\|+\|({\ce}^{-1})_{\infty}\|)^{\epsilon}.
$$
(For the definition of $(\cdots )_{\infty}$ and the norms see section 8.5.)
If we let $\l$ vary in a compact subset in the 
tube domain over this neighborhood
then there is a constant $c_2\ge c_1$ such that 
$$
c_1(1+\|\Im(\l)+\check{\eta}(m)\|+\|({\ce}^{-1})_{\infty}\|)^{\epsilon}
\le c_2(1+\|\check{\eta}(m)\|+\|({\ce}^{-1})_{\infty}\|)^{\epsilon}
$$
For $\Re(\l)$ sufficiently small and $\a\in \D_{P'}$ we have 
$$
\langle \Re(\l)+\l_L,\a\rangle\ge \frac{1}{2}\langle \l_L,\a\rangle >0.
$$
Therefore, there are 
$c_3,c_4>0$ such that for all such $\l$ and $m\in M_{\R}$
we have
$$
c_3\le 
\left|
 \prod_{\a\in \D_{P'}}
 \frac{\langle \l+\l_L-i\check{\eta}(m),\a\rangle}{
 \langle \l+\l_L-i\check{\eta}(m),\a\rangle+1}  
\right|  \le c_4
$$
Putting everything together, we can conclude that there is a neighborhood
$B$ of ${\bf 0}$ in  $ PL(\Sigma)_{\R}\oplus X^*(P)_{\R}$ such that 
for $(\p,\l) $ in a compact subset ${\bf K}$ 
of the tube domain over $B$ we have
$$
\left|
\prod_{e\in \Sigma_1}
L_f(\c_e, 1+(\p+\p_L+im)(e))
\prod_{e\in \Sigma_1-\Sigma_1'}
\frac{(\p+im)(e)}{(\p+im)(e)+1}\right|
$$
$$
\times
\left| E^G_P(\tl-i\check{\eta}(m),\ce^{-1})
\prod_{\a\in \D_P-\D_{P'}}
\frac{ \langle \l-i\check{\eta}(m),\a\rangle}{\langle 
\l-i\check{\eta}(m),\a\rangle+1}
\right|
\le c'({\bf K})(1+\|m+m_{\infty}(\c)\|)^{\epsilon},
$$
where $\|\cdot\|$ is a norm on $M_{\R,\infty}$.
On the other hand, by Lemma 3.5, we have
$$
\left|\zeta_{\Sigma}(\c,\tp +im)\right|\le
c''({\bf K})\prod_{v|\infty}
\left\{
\sum_{\dim \s =d}\,\,\prod_{e\in \s\cap \Sigma_1}
\frac{1}{(1+|e(m+m_v(\c))|)^{1+1/d}}
\right\}
$$
Now we may choose $\epsilon $ and $c_5$ such that
$$
(1+\|m+m_{\infty}(\c)\|)^{\epsilon}\le c_5\prod_{v|\infty}
\prod_{e\in \s_v\cap \Sigma_1}(1+|e(m+m_v(\c))|)^{1/2d}
$$
for any system $(\s_v)_{v|\infty}$ of $d$-dimensional cones.
This gives the claimed estimate.
\hfill $\Box$

\subsection*\no
{\bf 5.5}\hskip 0,5cm
To begin with, we let
$$
M'_{\R}:=\{ m\in M_{\R}\,\,|\,\, e(m)=0\, \forall e\in \SL\,\,
{\rm and}\,\, \langle \check{\e}(m),\a\rangle =0\,\, \forall\, \a\in 
 \D_P -\D_{P'}  \},
$$
$$
M'=M'_{\R}\cap M.
$$
Then $M'_{\R}=M'\otimes \R$ and $M''=M/M'$ is torsion free. Put
$d'={\rm rank}(M')$, $d''={\rm rank}(M'')$.

The connection with sections 5.2 and 5.3 is as follows:
$$
E=(PL(\Sigma)_{\R}\oplus X^*(P)_{\R})/M'_{\R},
$$
$$
V=M''\otimes \R, A=M'', \nu = d'',
$$
the set of linear forms $l_1,...,l_m$ is given as follows

\begin{equation}
\label{5.5.1-0}
(\p,\l)+M'_{\R}\mapsto \p(e), \,\, e\in \SL,
\end{equation}

\begin{equation}
(\p,\l)+M'_{\R}\mapsto \langle \l,\a\rangle, \a\in \DL
\label{5.5.1-1}
\end{equation}

The measure $dv=dm''$ on $V=M''_{\R}$ is normalized by $M''$, $dm=dm'dm''$,
where $dm$ (resp. $dm'$) is the Lebesgue measure on $M_{\R}$
(resp. $M'_{\R}$) normalized by $M$ (resp. $M'$).

Fix a convex open neighborhood of
${\bf 0}$ in $PL(\Sigma)_{\R}\oplus X^*(P)_{\R}$
for which Lemma 5.4 is valid. Denote by $B$ the image of this 
neighborhood in $E$. This is an open and 
convex neighborhood of ${\bf 0}$. 
Using Lemma 5.4 we see that
$$
g(\p,\l)=
\mu_{T'}\int_{M'_{\R}}\frac{1}{\kappa^{d''}}
\left\{\sum_{\chi\in {\cal U}_T}
\hat{H}_{\Sigma}(\chi,\tp+im')E^G_P(\tl-i\check{\eta}(m'),\ce^{-1}) 
h_L(\p,\l)\right\}dm'
$$
is a holomorphic function on $T_B$ (here $ \mu_{T'}=1/(2\pi \kappa)^{d'}$). 
(We use the invariance of $g$ under $iM_{\R}'$ and Cauchy-Riemann 
differential equations to check that $g$ is actually a function on 
$T_B$.)
Hence, 
$$
f(\p,\l):=g(\p,\l)h_L(\p,\l)^{-1} 
$$
is a meromorphic function on $T_B$. 

Let $E^{(0)}$ be the common kernel of all maps (\ref{5.5.1-0}, \ref{5.5.1-1}). 
Note that there is an exact sequence 
$$
0\ra M'_{\R}\ra E^{(0)}\ra \langle \L(L)\rangle \ra 0
$$
which implies 

\begin{equation}
b(L)={\rm codim}\,\,\L(L)=m-d'',
\label{5.5.2}
\end{equation}
where $m=\#(\SL ) +\#(\DL)$. 

By construction, $M_{\R}''\cap E^{(0)}=\{{\bf 0}\}$. Let $dy$ be the
Lebesgue measure on $E_0^{\vee}$ normalized by the lattice
generated by the linear forms (\ref{5.5.1-0}, \ref{5.5.1-1}).
Denote by $E^+_0$ the closed simplicial cone in $E_0$ defined by these linear
forms, and by
$ \pi_0\,:\, E\ra E_0$ the canonical projection. It is easily seen that
$\pi_0(M_{\R}'')\cap E^+_0=\{{\bf 0}\}$ (using the exact sequence above).
Let 
$$
\psi_0\,:\, E_0\ra P:=E_0/\pi_0(M_{\R}'')
$$
be the canonical projection and put
$$
\L=\psi_0(E^+_0),
$$
$$
B^+=B\cap \{ (\p,\l)\in E\,|\, \p(e)>0\,\,\forall \,e\in \SL, 
\langle \l,\a\rangle >0\,\,\forall\, \a\in \DL \}.
$$
By the following theorem the function $f\in {\cal M}(T_B)$ is 
distinguished with respect to $M_{\R}''$ and the set of linear forms
$(\ref{5.5.1-0}, \ref{5.5.1-1})$. Therefore, we can define
$\tilde{f}_{B^+}\,:\, T_{B^+}\ra \C$ by
$$
\tilde{f}_{B^+}(z)=\frac{1}{(2\pi)^{d''}}\int_{M_{\R}''}
f\left(z+i(m'',-\check{\eta}(m''))\right)dm''.
$$

\bigskip
\no
{\bf Theorem 5.5}
\hskip 0,5cm 
{\it
a) $f$ is a distinguished function with respect to $M_{\R}''$ and the set
of linear forms $(\ref{5.5.1-0}, \ref{5.5.1-1})$.

b) There exist an open neighborhood $\tilde{B}$ of ${\bf 0}$ 
containing $B^+$ and linear forms $\tilde{l}_1,...,\tilde{l}_{\tilde{m}}$
which vanish on $M_{\R}''$ such that 
$$
\tilde{f}_{B^+}(z)\prod_{j=1}^{\tilde{m}}\tilde{l}_j(z)
$$
has a holomorphic continuation to $T_{\tilde{B}}$ and $g({\bf 0})\neq 0$.
}
\bigskip

{\em Proof.}
a) Define $c_0\,:\, M_{\R,\infty}\ra \R_{\ge 0}$ by
$$
c_0((m_v)_v)=\prod_{v|\infty}
\{ \sum_{\dim \s =d }\,\,\prod_{e\in \s\cap \Sigma_1}
 \frac{1}{(1+|e(m_v)|)^{1+1/2d}}\}.
$$
Let ${\cal F}\subset M_{\R,\infty}^{1}$ be the cube spanned by
a basis of the image of ${\cal U}_T$ in  $M_{\R,\infty}^{1}$.
Let $c'>0$ such that for all $m_{\infty}(\chi)\in  M_{\R,\infty}^{1}$
($\chi\in {\cal U}_T)$ and all $m^1\in {\cal F}$
$$
c_0(m_{\infty}(\chi))\le c'c_0(m_{\infty}(\chi)+m^1).
$$
Let $dm^1$ be the Lebesgue measure on $M_{\R,\infty}^{1}$
normalized by the image of ${\cal U}_T$. Then for $m\in M_{\R}$
$$
\int_{M'_{\R}}\left\{\sum_{\c\in {\cal U}_T}c_0(m_{\infty}(\c)+m'+m)
\right\}dm'\le c(m \hskip 0,2cm {\rm mod} \hskip 0,2cm M'_{\R}),
$$
where $c\,:\, M^{''}_{\R}\ra \R_{\ge 0}$ is defined by
$$
c(m^{''}):={\rm cl}_F^dc'\int_{M'_{\R}}\int_{M_{\R,\infty}^{1}}
c_0(m'+m^1+m^{''})dm^1dm'.
$$

\no
By Lemma 5.4, for any compact subset ${\bf K}$ of $T_B$ there is 
a $c({\bf K})>0$  such that
$$
|g(z+im'')|\le c(m'')
$$
for all $z\in {\bf K}$ and $m''\in M_{\R}''$.
Obviously, $c$ can be integrated over any subspace $U$ of $M_{\R}''$.
It remains to show that for any $m''\in M_{\R}'' - U$ one has
$$
\lim_{\tau \ra \pm \infty}\int_Uc(\tau m''+u)du=0.
$$
This exercise will be left to the reader.

b) The first part concerning the meromorphic continuation and singularities
of $\tilde{f}_{B^+}$ is the content of Theorem 5.2. The relation 
$$
\lim_{s\ra 0}s^{b(L)}\tilde{f}_{B^+}(s\hat{L})=
g({\bf 0}){\cal X}_{\L}(\psi_0(\hat{L}))
$$
is satisfied by Theorem 5.3 and (\ref{5.5.2}). 
It will be shown in Section 6 that $g({\bf 0})\neq 0$. 
\hfill $\Box$

\subsection*\no
{\bf 5.6}\hskip 0,5cm 
The main theorem of our paper is:

\bigskip
\no
{\bf Theorem 5.6}\hskip 0,5cm
{\it
Let $L$ be a line bundle on $Y$ which lies in the interior of the cone
of effective divisors. Then there exists a metrization ${\cal L}$ of $L$ with 
the following properties:

a) The height zeta function
$$
Z_{Y^o}({\cal L},s)=\sum_{y\in Y^o(F)}H_{\cal L}(y)^{-s}
$$
is holomorphic for $\Re (s)>a(L)$ and it can be continued meromorphically
to a half-space $\Re(s)>a(L)-\d$ for some $\d>0$. In this half-space 
it has a pole of order $b(L)$ at $a(L)$ and no other poles.

b) For the counting function one has the following asymptotic
relation
$$
N_{Y^{o}}({\cal L}, H)= c({\cal L})H^{a(L)}(\log H)^{b(L)-1}(1+o(1))
$$
for $H\ra \infty$ with some constant $c({\cal L})>0$. 
}
\bigskip

{\em Proof.}
a) By construction, $\hat{L}=\frac{1}{a(L)}(\p_{\Sigma}+\p_L,2\l_P+\l_L)$ 
is mapped onto $[L]$ by $\psi$. Hence $Z_{Y^o}({\cal L},s)$ converges 
absolutely for $\Re(s)>a(L)$, where ${\cal L}$ is the metrization mentioned
in Proposition 4.4. By the same proposition, 
$$
Z_{Y^o}({\cal L}, s+a(L))=\mu_T
\int_{M_{\R}}\{\sum_{\chi\in {\cal U}_T}f_L(\chi,im)\}dm
$$
where 
$$
f_L(\chi,im):= 
$$
$$
\hat{H}_{\Sigma}(\chi,\frac{s}{a(L)}(\p_{\Sigma}+\p_L) +
\p_{\Sigma}+\p_L+im)
E^G_P(\frac{s}{a(L)}(2\rho_P+\l_L)+\rho_P+
\l_L-i\check{\eta}(m),\ce^{-1})
$$
for all $s$ with $\Re(s)>0$. However, this is just
$$
\frac{1}{(2\pi)^{d''}}
\int_{M_{\R}''}f\left(s\hat{L}+i(m'',-\check{\eta}(m''))\right)dm''=
\tilde{f}_{B^+}(s\hat{L})
$$
with $f,B^+$ and $\tilde{f}_{B^+}$ introduced in the preceding section.
By Theorem 5.5, $\tilde{f}_{B^+}$ extends to a meromorphic function on 
a tube domain over a neighborhood of ${\bf 0} $ and in this tube domain the 
only singularities are the hyperplanes defined over $\R$. 
Hence there is a $\d>0$ 
such that $Z_{Y^o}({\cal L},s+a(L))$ extends to a 
meromorphic function in the half-space $\Re(s)>-\d$ and the only 
possible pole is in $s=0$ and its order is exactly $b(L)$ 
(Theorems 5.3 and 5.5). 

b) This result follows from a Tauberian theorem 
(cf. \cite{De}, Th\'eor\`eme III or \cite{Sh}, Problem 14.1 (in the constant 
stated there the factor $\frac{1}{k_0}$ is missing)).
\hfill $\Box$

\section{Non-vanishing of asymptotic constants}
\label{6}

\subsection*\no 
{\bf 6.1} \hskip 0,5cm
This section is devoted to the proof of the non-vanishing of $g({\bf 0})$
claimed in Theorem 5.5.
All notations are as in sections 5.2-5.5.
The function $g(\p,\l)$ which has been defined in 5.5
is given by
$$
g(\p,\l)=\frac{\mu_{T'}}{\kappa^{d''}}\int_{M'}
\left\{\sum_{\chi\in {\cal U}_T}\hat{H}_{\Sigma}(\chi,\tp+im')
E^G_P\left(\tl -i\check{\eta}(m'),\ce^{-1}\right)h_L(\p,\l)
\right\}dm,
$$
where $\tp,\tl$ have been defined in 5.4.
The function $h_L(\p,\l)$ was defined in 5.4:
$$
h_L(\p,\l)=\prod_{e\in \SL}
\frac{\p(e)}{\p(e)+1}\prod_{\a\in \DL}
\frac{\langle \l,\a\rangle}{
\langle \l,\a\rangle +1}.
$$
The uniform convergency of the integral above in
any compact subset of $T_B$ (cf. Lemma 5.4) 
allows us to compute the limit 
$$
\lim_{(\p,\l)\ra {\bf 0}}\prod_{e\in \SL}\p(e)
\prod_{\a\in \DL}\langle \l,\a\rangle 
\hat{H}_{\Sigma}(\chi,\tp +im')
E^G_P\left( \tl -i\check{\eta}(m'),\ce^{-1}\right)
$$
first and then to integrate.
We shall show that this limit vanishes if there are
$e\in  \SL$ with $\chi_e\neq 1$ or
$\a\in  \DL$ with 
$\ce\circ\check{\a}\neq 1$.
Therefore, we may consider only $\c\in {\cal U}_T'$ where
$$
{\cal U}_T':=\{\c\in {\cal U}_T\,|\, \c_e=1\,\,\forall 
e\in  \SL,\,\,
\ce\circ\check{\a}=1 \,\forall\,\, \a\in
  \DL\}.
$$
Let $\eta'\,:\, P'\ra T$ be the uniquely defined homomorphism such that 
for all $\a\in \Delta_{P'}$ we have 
$\eta'\circ\check{\a}=\eta\circ\check{\a}$.

\bigskip
\no
{\bf Lemma 6.1}\hskip 0,5cm
{\it
\begin{equation}
\label{6.1.1}
g({\bf 0})
=\lim_{\p\ra {\bf 0},\,\p\in PL(\Sigma)^+ }\prod_{e\in \SL}\p(e)\cdot
\frac{\mu_{T'}}{\kappa^{d''}}\cdot
\frac{c_{P'}}{c_P} \int_{M'_{\R}}\{\sum_{\c\in {\cal U}_T'} 
\hat{H}_{\Sigma}(\c,\tp+im')
\end{equation}
$$
\hskip 8,5cm
\times 
E^G_{P'}\left(\l_L+\rho_{P'}-i\check{\e}'(m'),(\chi_{\eta'})^{-1}\right)\}dm'
$$
(cf. 8.4 for the definition of $c_{P'}$ and $c_P$). 
}
\bigskip

{\em Proof.}
Recall that (cf. (\ref{3.5.1}))
$$
\hat{H}_{\Sigma}(\chi,\tp+im')=\z_{\Sigma}(\chi, \tp +im')
\prod_{e\in \Sigma_1}L_f(\chi_e,
1+\p(e)+\p_L(e)+ie(m'))
$$
and that $\zeta_{\Sigma}(\chi,\tp +im')$ is regular for $\p$ in a tube domain
over a neighborhood of ${\bf 0}$ (Lemma 3.5).
For $e\in \Sigma_1'$ we have $\p_L(e)>0$, hence we see that the function
$$
L_f(\chi_e, 1+\p(e)+\p_L(e)+ie(m'))
$$ 
is holomorphic for 
$\Re (\p)$ in a neighborhood of
${\bf 0}$. 
Let $e\in \SL$. If $\chi_e\neq 1$ 
then  the restriction of
$\chi_e$ to ${{\bf G}_m({\bf A})^1}$ is non-trivial
(by our construction of the 
embedding ${\cal U}_T\ra {\cal A}_T$, cf. 3.3), hence
$$
\p\mapsto L_f(\chi_e, 1+\p(e))
$$
is an entire function and $\p(e)L_f(\chi_e, 1+\p(e))$ tends to $0$ as
$\p\ra {\bf 0}$. 

For $\a\in \Delta_{P'}$ we have $\langle \l_L,\a\rangle >0$, hence
$\l_L$ is contained in $ X^*(P')^+$.
Let $\a\in \DL$. If $\ce\circ \check{\a}
\neq 1$ then $\ce\circ \check{\a}$ 
restricted to ${{\bf G}_m({\bf A})^1}$ is non-trivial 
and therefore
$$
\prod_{\a\in  \DL}\langle \l,\a\rangle E^G_P\left(
\l+\rho_P+\l_L-i\check{\eta}(m'),\ce^{-1}\right)
$$
vanishes as $\l\ra {\bf 0}$ (cf. Proposition 8.3).
We have shown that it suffices to 
take the sum over all $\chi\in {\cal U}_T'$.
To complete the proof, note that for $\chi\in {\cal U}_T'$ we have
(cf. Proposition (8.4))
$$
\lim_{\l\ra {\bf 0}}
\prod_{\a\in  \DL}\langle \l,\a\rangle
E^G_P\left( \tl-i\check{\eta}(m'),\ce^{-1}\right )
= \frac{c_{P'}}{c_P}E^G_{P'}\left(\l_L+\rho_{P'}-i\check{\eta}'(m'), 
(\ce')^{-1}\right).
$$
\hfill $\Box$

\subsection*\no
{\bf 6.2}\hskip 0,5cm
By the absolute and uniform convergence of
$$
\int_{M'_{\R}}\left\{\sum_{\chi\in {\cal U}_T'}
\hat{H}_{\Sigma}(\chi,\tp +im')
\prod_{e\in \SL}\p(e)\right\}dm'
$$
(cf. Lemma 3.5 and Theorem 3.6) and the convergence
of 
$$
\sum_{\g \in P'(F)\ba G(F)}e^{\langle \l_L+2\rho_{P'},H_{P'}(\g)\rangle}
$$
we may change summation and integration in (\ref{6.1.1}) and get for  all 
$\p \in PL(\Sigma)^+$ 
$$
\frac{\mu_{T'}}{\kappa^{d''}}\cdot
\frac{c_{P'}}{c_P} \int_{M'_{\R}}\sum_{\c\in {\cal U}_T'}
\hat{H}_{\Sigma}(\c,\tp+im')
E^G_{P'}\left(
\l_L+\rho_{P'}-i\check{\eta}'(m')(\chi_{\eta'})^{-1}\right)dm'
$$
$$
=\frac{c_{P'}}{c_P\kappa^{d''}}\sum_{\g\in P'(F)\ba G(F)}
e^{\langle \l_L+2\rho_{P'},H_{P'}({\g}')\rangle}
\mu_{T'}
\int_{M'_{\R}} \sum_{\c\in \U_T^{'}}\hat{H}_{\Sigma}(\c,\tp+im')
(\c^{m'}\c)^{-1}(\e'(p_{\g}'))  dm'
$$
where $\g =p_{\g}'k_{\g}$ as above. 
Let $I$ be the image of the homomorphism 
$$
\prod_{e\in \SL}{\bf G}_m({\bf A})\times
\prod_{\a\in \DL}{\bf G}_m({\bf A})\ra T({\bf A})
$$
induced by $M\ra \Z^{\SL}\oplus 
\Z^{\DL}$, 
$$
m\mapsto
\left((e(m))_{e\in \SL},(\langle -m\circ\eta,\a\rangle)_{
\a\in  \DL})\right).
$$
Then $M_{\R}'\oplus  {\cal U}_T'$ is precisely the set of characters 
$T({\bf A})\ra S^1$ which are trivial on 
$T(F){\bf K}_TI$. Put $T'=\Spec (F[M'])$ and $ T''=\Spec (F[M'']).$
Then there is an exact sequence 
$$
1\ra T''\ra T\ra T'\ra 1.
$$
Note that $I\subset T''({\bf A})$. Denote by ${\bf K}_{T'}$ 
(resp. ${\bf K}_{T''}$) the maximal compact subgroup of $T'({\bf A})$
(resp. of $T''({\bf A})$).
The linear forms
$$
m\mapsto e(m),\,\,\, e\in \SL,
$$
$$
m\mapsto -\langle m\circ \eta,\a\rangle, \,\,\,
\a\in \DL,
$$
when considered as functions on $M''$,
generate a sublattice of finite index in
$N''=\Hom (M'',\Z)$. 
This shows that there is a $q>0$ such that the image of the $q$-th power
homomorphism $T''({\bf A})\ra T''({\bf A})$, 
$t\mapsto t^q$, is contained in $I$. 
If $v$ is any archimedean place of $F$ the connected component of one
in $T''(F_v)$ is therefore contained in $I$. Consequently, 
$$
T''(F)\cdot\prod_{v|\infty}T''(F_v)
\prod_{v\nmid\infty }{\bf K}_{T'',v}
\subset T''(F){\bf K}_{T''}\cdot I
$$
and the left hand side is of finite index in $T''({\bf A})$.
Put
$$
{\cal A}_T'=\{\c\in {\cal A}_T\,\,|\,\, \c=1 \,\, {\rm on}\,\,
T(F){\bf K}_TT''({\bf A})\}.
$$
We observe that
$$
{\cal A}_T'\simeq {\cal A}_{T'}=(T'({\bf A})/T'(F){\bf K}_{T'})^*\subset
M'_{\R}\oplus{\cal U}_T'
$$ 
We denote by 
$$
\Ind(L)=(M'_{\R}\oplus{\cal U}_T')/{\cal A}_T'
$$
and by $\ind(L)$  the  order of $\Ind(L)$.
Put ${\cal U}_{T'}={\cal U}_T\cap {\cal A}_T'$
(then $ {\cal A}_T'=M_{\R}'\oplus
{\cal U}_{T'}$).
Thus we may write
$$
\int_{M'_{\R}}\left\{\sum_{\c\in {\cal U}'_T}
\hat{H}_{\Sigma}(\c,\tp+im')(\c^{m'}\c)^{-1}(\e(p_{\g}'))
\right\}dm'
$$
$$
= \sum_{\c\in \Ind(L) }
\int_{M'_{\R}}\left\{\sum_{\c'\in {\cal U}_{T'}}
\hat{H}_{\Sigma}(\c'\c^{m'}\c,\tp)
(\c'\c^{m'}\c)^{-1}(\e'(p'_{\g}))\right\}dm'
$$
For $\c\in M'_{\R}\oplus {\cal U}_T'$ and $x\in T'({\bf A})$
we consider  the function
$$
x\mapsto \int_{T''({\bf A})}H_{\Sigma}(xt\e'(p_{\g}'), \tp)\c(t)dt
$$
(the Haar measure $dt$ on $T''({\bf A})$ 
is defined as $dx$ on $T({\bf A})$, cf. 3.3).
The same argument as in the proof of Proposition 3.4 shows that this 
function is absolutely integrable over $T'(F)$ if $\p\in PL(\Sigma)^+$. 
The  Fourier transform 
for $\c'\in {\cal A}_T'$
$$
\int_{T'({\bf A})}\left(\int_{T''({\bf A})}H_{\Sigma}(xt\eta'(p_{\g}'),\tp)
\c(t)dt\right)\c'(x)dx
$$
is absolutely integrable over ${\cal A}'_T$.
Using Poisson's summation formula twice we get
$$
\sum_{\c\in \Ind(L)}\mu_{T'}
\int_{M_{\R}'}\left\{\sum_{\chi'\in {\cal U}_{T'}}
\hat{H}_{\Sigma}(\c^{m'}\c'\c,\tp)
(\c'\c^{m'}\c)^{-1}(\e'(p'_{\g}))\right\}dm'
$$
$$
=\sum_{\c\in \Ind(L) }
\sum_{x\in T'(F)}\int_{T''(\A)}
H_{\Sigma}(xt\e'(p_{\g}',\tp)\c(t)dt=
\sum_{x\in T'(F)} \ind(L)
\int_{T''(F) {\bf K}_T I}H_{\Sigma}(xt\e'(p_{\g}'),\tp)dt.
$$

\no
Now we collect all the terms together. 

\bigskip
\no
{\bf Lemma 6.2}\hskip 0,5cm
{\it 
The constant $g({\bf 0})$ is equal to
$$
\frac{\ind(L) c_{P'}/c_{P}}{\kappa^{d''}} 
\sum_{\g\in P'(F)\ba G(F)}e^{\langle \l_L+2\rho_{P'},H_{P'}(\g)\rangle}
\times
$$
$$
\lim_{\p\ra {\bf 0},\,\p\in PL(\Sigma)^+} \prod_{e\in \SL}\p(e) 
\sum_{x\in T'(F)}\int_{T''(F){\bf K}_{T''}I}
H_{\Sigma}(x\e'(p_{\g}')t,\tp)dt.
$$
}
\bigskip

\no
{\bf Lemma 6.3}\hskip 0,5cm
{\it
The limit
$$
\lim_{\p\ra {\bf 0},\p\in PL(\Sigma)^+}
\prod_{e\in \SL}\p(e)\int_{T''({\bf A})}
H_{\Sigma}(t,\tp)dt
$$
exists and is positive.
}

\bigskip
{\em Proof.}
Consider the embedding $N''_{\R}\ra N_{\R}$ and let
$$
\Sigma'':=\{\s\cap N_{\R}''\,|\, \s\in \Sigma\}.
$$
This is a complete fan in $N_{\R}''$ which consists of rational polyhedral 
cones, but which is not necessary a regular fan. We can obtain a regular fan
by subdivision of the cones into regular ones (cf. 
\cite{KKMS-D}, ch. I, \S 2, 
Theorem 11). This gives us a complete regular 
fan $\tilde{\Sigma}''$ such that any 
cone in $\tilde{\Sigma}''$ is contained in a cone of $\Sigma''$. Denote by 
$\tilde{\Sigma}_1''$ the set of primitive integral generators of 
the one-dimensional cones in $\tilde{\Sigma}''$. 
Computing the integral as in section 3.4 we get
$$
\int_{T''({\bf A})}H_{\Sigma}(t,\tp )dt=\z_{\tilde{\Sigma}''}(1,\tp  )
\prod_{\tilde{e}\in \tilde{\Sigma}_1''}L_f(1,\tp(\tilde{e}))
$$
(cf. (\ref{3.5.1})), where
$\zeta_{\tilde{\Sigma}''}(1,\tp)$ is regular in a 
neighborhood of $\p={\bf 0}$ and positive for $\p={\bf 0}$. 
Let $\tilde{e}\in N$ and $\s\in \Sigma'$ 
be a cone containing $\tilde{e}$. Write
$\tilde{e}=\sum_{e\in \s\cap \Sigma_1}t_e\cdot e$ ($ t_e\in \Z_{\ge 0}$).
Suppose
$$
1=(\p_{\Sigma}+\p_L)(\tilde{e})=\sum_{e\in \s\cap \Sigma_1}t_e(1+\p_L(e)).
$$
Then $t_e=0$ for all $e\in \s\cap \Sigma_1'$ 
because $\p_L\in \sum_{e\in \Sigma_1'}
\R_{>0}\p_e$ (cf. 5.4).
Hence, 
$$
\tilde{e}=\sum_{e\in \s\cap \SL}t_e\cdot e
$$ 
and
$\p_{\Sigma}(\tilde{e})=1$ implies $\tilde{e}\in \SL$.
Therefore,
$$
\lim_{\p\ra {\bf 0},\,\p\in PL(\Sigma)^+}\prod_{e\in \SL}\p(e)
\prod_{\tilde{e}\in \tilde{\Sigma}_1^{''}}L_f(1, \tp(\tilde{e}))
$$
$$
= \{\prod_{e\in \SL}\lim_{\p\ra{\bf 0},\,\p\in PL(\Sigma)^+}\p(e)
L_f(1,1+\p(e))\}
\prod_{\tilde{e}\in \tilde{\Sigma}_1''- (\SL) }
L_f\left(1,(\p_{\Sigma}+\p_L)(\tilde{e})\right)
$$
and this is a positive real number. 
\hfill $\Box$

\bigskip 
\no
In theorem 5.5 we claimed the non-vanishing of $g({\bf 0})$. 
We are now in the position to prove

\bigskip
\no
{\bf Corollary 6.4}\hskip 0,5cm
{\it
$$
g({\bf 0})>0.
$$
}
\bigskip

{\em Proof.}
By Lemma 6.2 it is enough to show that
$$
\lim_{\p\ra {\bf 0}, \p\in PL(\Sigma)^+}\prod_{e\in \SL}\p(e)
\int_{T''(F){\bf K}_{T''}\cdot I}H_{\Sigma}(t,\tp)dt
$$
is positive. Let $t_1,...,t_{\nu}\in T({\bf A})$ such that
$$
T''({\bf A})=\bigcup_{j=1}^{\nu}t_jT''(F){\bf K}_{T''}\cdot I.
$$
Then there exists a constant $c>0$ such that for all $t\in T''({\bf A})$
and
$j=1,...,\nu$
we have
$$
H_{\Sigma}(tt_j,\tp)\le \frac{c}{\nu}H_{\Sigma}(t,\tp).
$$
Hence we can estimate
$$
\int_{T''({\bf A})}H_{\Sigma}(t,\tp)dt=\sum_{j=1}^{\nu}
\int_{T''(F){\bf K}_{T''} I}H_{\Sigma}(tt_j,\tp)dt
$$
$$
\le c\int_{T''(F){\bf K}_{T''} I}H_{\Sigma}(t,\tp)dt.
$$
Lemma 6.3 allows us to conclude that the limit above is indeed positive.
\hfill $\Box$

\section{Technical theorems}
\label{7}

\subsection*\no
{\bf 7.1}\hskip 0,5cm
Let $(A, V, \L) $ be a triple consisting of
a free abelian group
$A$ of rank $d$, a $d$-dimensional real vector space
$V:= A \otimes {\bf R}$ containing $A$ as a sublattice of
maximal rank, and  a closed strongly convex
polyhedral $d$-dimensional cone
$\Lambda \subset A_{\R}$ such that $\Lambda \cap - \L =\{{\bf  0}\}$. 
Denote by  $\L^{\circ}$ the interior  of $\L$. 
Let $( A^{\vee} , V^{\vee }, \L^{\vee }) $ be the triple
consisting of the dual abelian group
$A^{\vee } = {\rm Hom}(A, \Z)$, the dual real vector space
$V^{\vee } = {\rm Hom}(V, \R)$ and the  dual cone
$\L^{\vee } \subset V^{\vee }$.
We normalize the Haar measure $ dy$ on $V^{\vee }$
by the condition:
${\rm Vol}(V^{\vee }/A^{\vee })=1$.

We denote by $\chi_{\L}(v)$ the set-theoretic characteristic
function of the cone $\L$  and by ${\cal X}_{\L}(v)$ the 
Laplace transform of the set-theoretic characteristic function 
of the dual cone
$$
{\cal X}_{\L}(v) =\int_{V^{\vee }}\chi_{\L^{\vee}}(y)e^{-\langle v,y
\rangle } dy=
\int_{{\L}^{\vee }} e^{- \langle v, y
 \rangle}  dy, 
$$
where ${\Re }(v) \in {\L}^{\circ}$ (for these $v$ 
the integral converges absolutely).

Consider a complete regular fan $\Sigma$ on 
$V$, that is, a subdivision of the real space
$V$ into a finite set of convex rational simplicial cones,
satisfying certain conditions 
(see \cite{BaTschi1}, 1.2). 
Denote by $\Sigma_1$ the set of primitive generators
of one dimensional cones in $\Sigma$. 
Denote by $PL(\Sigma)_{\R}$ the vector space of 
real valued piecewise linear functions on $V$ and by
$PL(\Sigma)_{\C}$ its complexification. 

\begin{prop}(\cite{BaTschi1}, Prop. 2.3.2, p. 614)
For any compact ${\bf K}\subset PL(\Sigma)_{\C}$ with the
property that $\Re(\varphi(v))>0$ for all $\varphi \in {\bf K}$
and $v\neq {\bf 0}$ there exists a constant $\kappa({\bf K})$
such that
for all $\varphi\in {\bf K}$ and all $y\in V^{\vee}$ 
the following inequality holds:
$$
\left|\int_{V}e^{-\varphi(v) -i \langle v,y \rangle}dv\right|
\le
\kappa({\bf K})\sum_{\dim \s=d}\frac{1}{\prod_{e\in \s}(1+
|\langle e,y \rangle |)^{1+1/d}}.
$$
\label{estimate-fan}
\end{prop}

\bigskip\no
{\bf 7.2}\hskip 0,5cm
Let $H\subset V$ be a hyperplane with $H\cap \L=\{{\bf 0}\}$.
Let $y_0\in V^{{\vee }}$ with 
$H=\Ker (y_0)$, such that for all $x\in \L\,:\, y_0(x)\ge 0$.
Then $y_0$ is in the interior of $\L^{\vee}\subset V^{\vee}$. 
Let $x_0\in \L^{\circ}$ and let 
$$
H'=\{y'\in V^{\vee }\,|\, y'(x_0)=0\}.
$$
We have $V^{\vee}= H'\oplus \R y_0$. Define $\p\,:\, H'\ra \R$ by
$$
\p(y')=\min \{ t\,|\, y'+ty_0\in \L^{{\vee }}\}.
$$
The function $\p$ is piecewise linear with respect to a
complete fan of $H'$. Taking a subdivision, if necessary,
we may assume it to be regular. 

\begin{prop}
The function ${\cal X}_{\L}(u)$ is 
absolutely integrable over any linear subspace  $U\subset H$. 
\label{estimate-chi}
\end{prop}
 
{\em Proof.}
For $h\in H$ we have
$$
{\cal X}_{\L}(x_0+ih)=\int_{\L^{\vee}}e^{-y(x_0+ih)}dy=
\int_{H'}\int_{\p(y')}^{\infty}e^{-(y'+ty_0)(x_0+ih)}dtdy'
$$
$$
= \int_{H'}e^{-\p(y')}e^{-iy'(h)}dy'
$$
Therefore, $h\mapsto {\cal X}_{\L}(x_0+ih)$ is the Fourier transform
of the function
$y'\mapsto e^{-\p(y')}$ on $H'\simeq H^{\vee}$. 
The statement follows now from 7.1.
\hfill $\Box$

\bigskip
\no
{\bf 7.3}\hskip 0,5cm
The rest of this section is devoted to the proof of the meromorphic
continuation of certain functions which are holomorphic
in tube domains over convex finitely generated polyhedral
cones. In section \ref{5} we have already introduced 
the terminology and explained how this technical
theorem is applied to height zeta functions.

Let $E$ be a finite dimensional vector space over $\R$
and $E_{\C}$ its complexification. 
Let $V\subset E$ be a subspace. 
We will call a function  $c\,:\, V\ra \R_{\ge 0}$ 
sufficient if it satisfies the following conditions: 

(i) For any subspace $U\subset V$ and any $v\in V$ the function
$U\ra \R$ defined by $u\ra c(v+u)$ is measurable on $U$ and the 
integral 
$$
c_U(v):= \int_U c(v+u)du
$$
is always finite ($du$ is a Lebesgue measure on $U$). 

(ii) For any subspace $U\subset V$ and every $v\in V - U$ we have 
$$
\lim_{\tau \ra \pm \infty} c_U(\tau \cdot v)=0.
$$

Let $l_1,...,l_m\in E^{\vee}=\Hom_{\R}(E,\R)$ be linearly independent
linear forms. Put $H_j={\Ker } (l_j)$ for $j=1,..., m$.
Let $B\subset E$ be an open and convex neighborhood of ${\bf 0}$,
such that for all $x\in B$ and  all 
$j=1,...,m$ we have $l_j(x)>-1$. We denote by 
$T_B:=B+iE \subset E_{\C}$ the complex tube domain over $B$.
We denote by ${\cal M}(T_B)$ the set of meromorphic functions on 
$T_B$.

A meromorphic function $f\in {\cal M}(T_B)$
will be called {\it distinguished }  with respect to the data
$(V; l_1,...,l_m)$ if it satisfies the following
conditions:

(i) The function 
$$
g(z):=f(z)\prod_{j=1}^m\frac{l_j(z)}{l_j(z)+1}
$$
is holomorphic in $T_B$. 

(ii) There exists a sufficient function
$c\,:\, V\ra \R_{\ge 0}$ 
such that for any compact ${\bf K}\subset T_B$ there is a constant
$\kappa({\bf K})\ge 0$ such that for all $z\in {\bf K}$ and all
$v \in V$ we have
$$
|g(z+iv)|\le \kappa({\bf K})c(v).
$$

Let $C$ be a connected component of $B - \bigcup_{j=1}^m H_j$
and $T_C$ a tube domain over $C$.
We will consider the following integral:
$$
\tilde{f}_C(z):=\frac{1}{(2\pi )^d}\int_Vf(z+iv)dv.
$$
Here we denoted by
$d=\dim V$ and by $dv$ a fixed Lebesgue measure on $V$. 

\begin{prop}
Assume that $f$ is an distinguished function 
with respect to $(V;l_1,...,l_m)$. Then the following holds: 

a) $\tilde{f}_C\,:\, T_{C}\ra \C$ is a holomorphic function.

b) There exist an open and convex neighborhood  $\tilde{B}$
of ${\bf 0}$, containing $C$, and linear forms 
$\tilde{l}_1,...,\tilde{l}_{\tilde{m}}$, which vanish on $V$,
such that
$$
z\ra \tilde{f}_C(z)\prod_{j=1}^{\tilde{m}}\tilde{l}_j(z)
$$
has a holomorphic continuation to $T_{\tilde{B}}$. 
\label{cont}
\end{prop}

{\em Proof.}
a) Let ${\bf K}\subset T_C$ be a compact subset and let 
$\kappa({\bf K})\ge 0$ be a real number such that for all
$z\in {\bf K} $ and all $v\in V$ we have
$|g(z+iv)|\le \kappa({\bf K})c(v)$. 
Since ${\bf K}$ is a compact and $C$ doesn't intersect 
any of the hyperplanes $H_j$ there exist
real numbers $c_j\ge 0$ for any $j=1,...,m$, 
such that for all $z\in {\bf K}$ and $v\in V$ 
the following inequalities hold
$$
\left|\frac{l_j(z+iv)+1}{l_j(z+iv)}\right|\le c_j.
$$
Therefore, for $z\in {\bf K}$ and $v\in V$ we have
$$
|f(z+iv)|\le c_1\cdots c_m\kappa({\bf K})c(v).
$$
It follows that on every
compact ${\bf K}\subset T_C$ the integral
converges absolutely and uniformly to 
a holomorphic function $\tilde{f}_C$.

\medskip 
b) The proof proceeds by induction on $d=\dim V$.
For $d=0$ there is nothing to prove. Assume that $d\ge 1$ 
and let $v_0\in V - \{{\bf 0}\}$ be a vector such that both
$v_0,-v_0\in B$.
We define 
$B_1\subset B$ as the set of all vectors $x\in B$ which
satisfy the following two conditions:
the vector $x\pm v_0\in B $ and 
$|\frac{l_j(x)}{l_j(v_0)}|\le\frac{1}{2}$
for all $j\in \{1,...,m\}$ with $l_j(v_0)\neq 0$. 
The set $B_1$ is a convex open neighborhood of ${\bf 0}\in E$.
Fix a vector $x_0\in C$. Without loss of generality we
can assume that 
$$
\{1,...,m_0\}:=\{j\in \{1,...,m\}\,|\, l_j(v_0)l_j(x_0)<0\}
$$
with $0\le m_0\le m$.
For $j\in \{1,...,m_0\}, k\in \{1,...,{\widehat{j} },...,m\}$
we define
$$
l_{j,k}(x):=
l_k(x)-l_j(x)\frac{l_k(v_0)}{l_j(v_0)},\hskip 1cm 
H_{j,k}:={\Ker} (l_{j,k})\subset E.
$$
For all $j\in \{1,...,m_0\}$ we have that 
$(l_{j,k})_{1\le k\le m, k\neq j}$ is a set of 
linearly independent linear forms on $E$. 
Moreover, for all $x\in B_1$ and $j\in \{1,...,m_0\}$
we have 
$$
x-\frac{l_j(x)}{l_j(v_0)}v_0 =
\left(1-\frac{l_j(x)}{l_j(v_0)}\right)x+\frac{l_j(x)}{l_j(v_0)}(x-v_0)\in B,
$$
in the case 
that $\frac{l_j(x)}{l_j(v_0)}\ge 0$
and,  similarly, in the case that $\frac{l_j(x)}{l_j(v_0)}< 0$ 
$$
x-\frac{l_j(x)}{l_j(v_0)}v_0 =\left(1+\frac{l_j(x)}{l_j(v_0)}\right)x
+ \left(-\frac{l_j(x)}{l_j(v_0)}\right)(x+v_0)\in B.
$$
Therefore, for all $x\in B_1$ and $j\in \{1,...,m_0\},
k\in \{1,...,{\hat{ j}},...,m\}$ we have
$$
l_{j,k}(x)=l_k\left(x-\frac{l_j(x)}{l_j(v_0)}v_0\right)>-1.
$$
Let $C_1$ be a connected component of 
$$
B_1 -\left(
\bigcup_{j=1}^{m_0}\bigcup_{1\le k\le m, k\neq j } H_{j,k}
\cup \bigcup_{j=1}^mH_j\right),
$$
which is contained in $C$. 
For $z\in T_{C}$ we define
$$
h_C(z):=\frac{1}{2\pi}\int_{-\infty}^{+\infty}f(z+i\tau v_0)d\tau =
\frac{1}{2\pi i}\int_{{\rm Re}(\l)=0}f(z+\l v_0)d\l.
$$
As in (i) one shows that $h_{C}$ is a holomorphic function on $T_{C}$.
For $x\in B_1$ and $\l\in [0,1]$ we have
$$
x+\l v_0=(1-\l)x+\l(x+v_0)\in B.
$$
If for some $z = x+iy \in T_{C_1}$
($x\in C_1$) and $\l\in [0,1]+i\R, j\in \{1,...,m\}$
we have
$$
l_j(z+\l v_0)=0,
$$
then it follows that $l_j(x)+{\rm Re}(\l)l_j(v_0)=0$,
and therefore, $l_j(x)l_j(v_0)<0$
(since $l_j(x)$ has the same sign as $l_j(x_0)$). 
Consequently, $j\in \{1,...,m_0\}$.

For $z\in T_{C_1}$ and $j\in \{1,...,m_0\}$ we put
$$
\l_j(z):=-\frac{l_j(z)}{l_j(v_0)}.
$$
By our assumptions, we have $ 0<{\rm Re}(\l_j(z))<\frac{1}{2}$.
>From $\l_j(z)=\l_{j'}(z),$ with $j,j'\in \{1,...,m_0\}$ and
$j\neq j'$ it follows now that 
$$
l_{j'}\left(z-\frac{l_j(z)}{l_j(v_0)}v_0\right)=0.
$$
In particular, we have $l_{j,j'}({\rm Re}(z))=0$. This is not
possible, because $z\in T_{C_1}$.

Assume now that $x\in B_1$. We have, assuming that $l_k(v_0)\neq 0$,
that
$$
|l_k(x+v_0)|=|l_k(v_0)|\cdot \left|\frac{l_k(x)}{l_k(v_0)}+1\right|\ge
 |l_k(v_0)|\cdot |-\frac{1}{2}+1|=\frac{1}{2}|l_k(v_0)|.
$$
If $l_k(v_0)=0$ then we have $l_k(x+v_0)=l_k(x)$.

Fix a $z\in T_{C_1}$. Then there exist 
some numbers $c_1(z)\ge 0, c_2(z)\ge 0$
such that for all $\l\in [0,1], \tau \in \R$ with $|\tau |\ge c_1(z)$
we have
$$
|f(z+\l v_0+i\tau v_0)|
$$
$$
=|g(z+\l v_0+i\tau v_0)|\cdot |\prod_{j=1}^{m}(1+\frac{1}{
l_j(z+\l v_0 +i\tau v_0)})|\le c_2(z)\kappa (z+[0,1]v_0)c(\tau v_0).
$$
For any $z\in T_{C_1}$ we have therefore
$$
h_C(z)=-\sum_{j=1}^{m_0}{\rm Res}_{\l=\l_j(z)}(\l\ra f(z+\l v_0))
+ \frac{1}{2\pi i}\int_{\Re (\l)=1}f(z+\l v_0)d\l .
$$
Since $\l_j(z)\neq \l_{j'}(z)$ for $j,j'\in \{1,...,m_0\}$,
$j\neq j'$, $z\in T_{C_1}$, we have
$$
{\rm Res}_{\l=\l_j(z)}(\l\ra f(z+\l v_0))
$$
$$
= \lim_{\l\ra \l_j(z)} (\l-\l_j(z)) \prod_{k=1}^m
\frac{ l_k(z+\l v_0)+1}{l_k(z+\l v_0)} g(z+\l v_0)
$$
$$
= \frac{1}{l_j(v_0)}g\left(z-\frac{ l_j(z)}{l_j(v_0)}v_0\right)
\prod_{1\le k\le m, k\neq j} \frac{l_{j,k}(z)+1}{l_{j,k}(z)}.
$$
Put now for $j\in \{1,...,m_0\}, z\in T_{B_1}$
$$
f_j(z):= -\frac{1}{l_j(v_0)}g\left(z-\frac{l_j(z)}{l_j(v_0)}v_0\right)\cdot
\prod_{1\le k\le m,\, k\neq j}\frac{l_{j,k}(z)+1}{l_{j,k}(z)}
$$
and 
$$
g_j(z):= -\frac{1}{l_j(v_0)}\cdot g\left(z-\frac{l_j(z)}{l_j(v_0)}v_0\right).
$$

Let $V_1\subset V$ be a hyperplane which does not contain $v_0$. 
We want to show that the function $f_j(z)$ is distinguished with respect to
$(V_1;(l_{j,k})_{1\le k\le m,\, k\neq j})$. The function $f_j$ is 
meromorphic on $T_{B_1}$. Also, for all $x\in B_1$ and all 
$k\in \{1,...,{\hat j},...,m\}$ we have $l_{j,k}(x)>-1$. 
Further, we have that
$$
g_j(z)=f_j(z)\cdot \prod_{1\le k\le m,\,
 k\neq j}\frac{l_{j,k}(z)}{l_{j,k}(z)+1}
$$ is a holomorphic function on $T_{B_1}$.

Let ${\bf K}_1\subset T_{B_1}$ ($\subset T_B$) be a compact, and let
$$
{\bf K}(j):=\{z-\frac{l_j(z)}{l_j(v_0)}v_0\,|\, z\in {\bf K}_1 \}.
$$
This is a compact subset in $T_B$. Put 
$$
\kappa_j({\bf K}_1)=\frac{1}{|l_j(v_0)|}\kappa({\bf K}(j)),
$$
where $\kappa({\bf K}(j)) $ is a constant such that 
$|g(z+iv)|\le \kappa ({\bf K}(j))c(v)$ for all $z\in {\bf K}(j)$ and
all $v\in V$. For $v_1\in V_1$ we put
$$
c_j(v_1)=c\left(v_1-\frac{l_j(v_1)}{l_j(v_0)}v_0\right).
$$
Then for all $z\in {\bf K}_1$ and $v\in V_1$ we have
$$
|g_j(z+iv_1)|\le \kappa_j({\bf K}_1)c_j(v_1).
$$
Moreover, for any subspace $U_1\subset V_1$ and all $v_1\in V_1$ the
function $U_1\ra \R, \,\, u_1\ra c_j(u_1+v_1)$ is measurable and 
we have
$$
c_{j,U_1}(v_1):=\int_{U_1}c_j(v_1+u_1)du_1\,  < \infty .
$$
For all $v_1\in V_1 - U_1$ we have 
$$
\lim_{\tau\ra\pm \infty} c_{j,U_1}(\tau v_1)=0.
$$
This shows that $f_j$ is distinguished with respect to 
$(V_1; (l_{j,k})_{1\le k\le m,\, k\neq j} )$. 

For $z\in T_{B_1}$ we put
$$
f_0(z):=\frac{1}{2\pi i} \int_{\Re (\l)=1} f(z+\l v_0)d\l .
$$
If $l_k(v_0)\neq 0$ we have (as above) for all $x\in B_1$
the following inequality 
$$
|l_k(x+v_0)|\ge \frac{1}{2} |l_k(v_0)|.
$$
Therefore, we conclude that the function 
$$
g_0(z):= f_0(z)\prod_{1\le k\le m,\, l_k(v_0)= 0}\frac{l_k(z)}{l_k(z)+1}
$$
$$
=
\frac{1}{2\pi }\int_{-\infty}^{+\infty}
\{\prod_{1\le k\le m, l_k(v_0)\neq 0}
\frac{l_k(z+v_0+i\tau v_0)+1}{l_k(z+v_0+i\tau v_0)} \}
g(z+v_0+i\tau v_0)d\tau
$$
is holomorphic in $T_{B_1}$.

Further, we have for $z\in {\bf K}_1 $ (with ${\bf K}_1\subset T_{B_1}$
a compact) and $v_1\in V_1$ the inequality
$$
|g_0(z+iv_1)|\le \kappa_0({\bf K}_1)c_0(v_1),
$$
where $\kappa_0({\bf K}_1)$ is some suitable constant and 
$$
c_0(v_1):=\int_{-\infty}^{+\infty}c(v_1+\tau v_0)d\tau.
$$
Again, for any subspace $U_1\subset V_1$ and any $u_1\in V_1$ we have
that the map $U_1\ra \R$ given by 
$u_1\ra c_0(v_1+u_1)$ is measurable, and 
that 
$$
c_{0,U_1}(v_1)=\int_{U_1}c_0(v_1+u_1)du_1 \,<\infty.
$$
For all $v_1\in V_1 - U_1$ we have
$$
\lim_{\tau \ra \pm \infty} c_{0,U_1}(\tau v_1) =0.
$$
Therefore, $f_0$ is distinguished with respect to
$(V_1; (l_k)_{1\le k\le m, \, l_k(v_0)\neq 0})$.

The Cauchy-Riemann equations imply that $g_0$ is invariant
under $\C v_0$, that is, for all $z_1,z_2\in T_{B_1}$ with
$z_1-z_2\in \C v_0$ we have $g_0(z_1)=g_0(z_2)$. 
We see that $f_0$ is also invariant under $\C v_0$ (in this sense),
as well as $f_1,...,f_{m_0}$ (this can be seen from the explicit
representation of these functions).
For $z\in T_{C_1}$ we have
$$
h_C(z)=f_0(z)+\sum_{j=1}^{m_0} f_j(z)=\sum_{j=0}^{m_0}f_j(z).
$$
Moreover,  for such $z$ we have
$$
\tilde{f}_C(z)=\frac{1}{(2\pi)^{d-1}}\int_{V_1}h_C(z+iv_1)dv_1
=\sum_{j=0}^{m_0}\frac{1}{(2\pi)^{d-1}}\int_{V_1}f_j(z+iv_1)dv_1,
$$
where $dv_1d\tau =dv $ (and ${\rm Vol}_{d\tau}(
\{\l v_0\,|\l\in [0,1]\} )=1)$.

By our induction hypothesis, there exists an open and convex neighborhood
$B'$ of ${\bf 0}$ in $E$ and linear forms 
$\tilde{l}_1,...,\tilde{l}_{\tilde m}$, which vanish on $V$, such that
$$
\tilde{f}_C(z)\prod_{j=1}^{\tilde m} \tilde{l}_j(z)
$$
has a holomorphic continuation to $T_{B'}$. 
(Strictly speaking, from the induction hypothesis it follows 
only that the linear forms $\tilde{l}_1,...,\tilde{l}_{\tilde{m}}\in E^{\vee }$
vanish on $V_1$. But since the functions $f_0,...,f_{m_0}$
``live'' already on a tube domain in $(E/\R v_0)_{\C}$, it follows that
the linear forms are also $\R v_0$-invariant.)

Now we notice that $\tilde{f}_C(z)\prod_{j=1}^{\tilde m} \tilde{l}_j(z)$
is holomorphic on $T_C$. Let $\tilde{B}$ be the convex hull of $B'\cup C$. 
Then we have that 
$$
\tilde{f}_C(z)\prod_{j=1}^{\tilde m} \tilde{l}_j(z)
$$
is holomorphic on $T_{\tilde B}$  (cf. \cite{Hor}, Theorem 2.5.10).
\hfill $\Box$

\bigskip \no
{\bf 7.4}\hskip 0,5cm
Let $E,V,l_1,...,l_m$ and $ 
B\subset \{x\in E\,|\, l_j(x)>-1\hskip 0,2cm \forall 
\hskip 0,2cm j=1,...,m\}$ be as above.
Let $C$ be a connected component of $B - \bigcup_{j=1}^m H_j$.
Let $f\in {\cal M}(T_B)$ be an distinguished function with respect
to $(V; l_1,...,l_m)$. 

\begin{prop}
There exist an open convex neighborhood $\tilde B$ of ${\bf 0}$ in $E$
containing $C$ and linear forms $\tilde{l}_1,...,\tilde{l}_{\tilde{m}}\in
E^{\vee }$ vanishing on $V$ such that

a) for all $j\in \{1,...,\tilde{m} \} $ we have 
$ {\rm Ker}(\tilde{l}_j)\cap C=\emptyset$

b) $\tilde{f}_C(z)\prod_{j=1}^{\tilde{m}} \tilde{l}_j(z)$
has a holomorphic
continuation to $T_{\tilde{B}}$.
\end{prop}

{\em Proof.} By the proposition above, there exist linear forms
$\tilde{l}_1,...,\tilde{l}_{\tilde{m}}$ such that $
V\subset \cap_{j=1}^{\tilde{m}} {\rm Ker}(\tilde{l}_j)$ and
$\tilde{f}_C(z)\prod_{j=1}^{\tilde{m}}\tilde{l}_j(z)$ 
has a  holomorphic continuation to a tube domain $T_{\tilde{B}}$
over a  convex open neighborhood $\tilde{B}$ of ${\bf 0}\in E$
containing $C$. 
Suppose that there exist an $x_0\in C$ and a 
$j_0\in \{1,...,\tilde{m}\}$ such that $\tilde{l}_{j_0}(x_0)=0$.
Then the function
$$
\tilde{f}_C(z)\prod_{j\neq j_0}\tilde{l}_j(z)
$$
is still holomorphic in $T_{\tilde{B}'}$
with $\tilde{B}'= (B- {\rm Ker}(\tilde{l}_{j_0}))\cup C$. 
It is easy to see that $\tilde{B}'$ is connected.
The convex hull of $\tilde{B}'$ is equal to $\tilde{B}$. 
Therefore, already the function
$$
\tilde{f}_C(z)\prod_{j\neq j_0}\tilde{l}_j(z)
$$
is holomorphic on $T_{\tilde{B}}$ (cf. \cite{Hor}, loc. cit.).

\subsection*\no
{\bf 7.5}\hskip 0,5cm
As above, let $E$ be a finite dimensional  vector space over 
$\R$ and
let $l_1,...,l_m$ be linearly independent linear forms on $E$.
Put $H_j:={\rm Ker}(l_j)$, for $j=1,...,m$,
$E^{(0)}=\bigcap_{j=1}^m H_j$   and 
$E_{0}=E/E^{(0)}$. Let
$\pi_0\,:\, E\ra E_0$ be the canonical projection. Let $V\subset E$
be a subspace with $V\cap E^{(0)}=\{ {\bf 0}\}$, such that
$\pi_0 |_V\,:\, V\ra E_0$ is an injective map.
Let
$$
E^+_0:=\{ x\in E_0 \,|\, l_j(x)\ge 0, \hskip 0,1cm j=1,...,m\}
$$
and let $\psi\,:\, E_0\ra P := E_0/\pi_0 (V)$ be the canonical
projection.
We want to assume that $\pi_0(V)\cap E^+_0=\{{\bf 0}\}$, so that
$\L:=\psi (E^+_0)$ is a strictly convex polyhedral cone.

Let $dy$ be the Haar measure on $E_0^{\vee }$ normalized by
${\rm Vol}_{dy}(E^{\vee }_0/\oplus_{j=1}^m\Z l_j)=1$. 
Let $A\subset V$ be a lattice, and let $dv$ be a measure
on $V$ normalized by
${\rm Vol}_{dv}(V/A)=1$. 
On $V^{\vee }$ we have a measure $dy'$ normalized by $A^{\vee }$ and
a section of the projection $E_0^{\vee }\ra V^{\vee }$ gives a measure
$dy''$ on $P^{\vee }$ with $dy=dy'dy''$.

Let $B\subset E$ be an open and convex neighborhood of ${\bf 0}$,
such that for all $x\in B$ and $j\in \{1,...,m\}$ we have 
$l_j(x)>-1$.
Let $f\in {\cal M}(T_B)$ be a meromorphic function
in the tube domain over $B$ which is
distinguished with respect to $(V;l_1,...,l_m)$. 
Put 
$$
B^+=B\cap \{x\in E\,|\, l_j(x)>0, \hskip 0,1cm j=1,...,m \},
$$
$$
\tilde{f}_{B^+} (z)=\frac{1}{(2\pi)^d}\int_Vf(z+iv)dv
$$
where $d=\dim V$.

By 7.3, 
the function $\tilde{f}_{B^+}\,:\, T_{B^+}\ra \C$ is holomorphic
and it has a meromorphic continuation to a neighborhood
of ${\bf 0}\in E_{\C}$.
Put 
$$ 
g(z)=f(z)\prod_{j=1}^m\frac{l_j(z)}{l_j(z)+1}.
$$

\begin{prop}
For $ x_0\in B^+ $ we have
$$
\lim_{s\ra 0}s^{m-d}\tilde{f}_{B^+}(sx_0)=g({\bf 0}){\cal X}(\psi (x_0)).
$$
\end{prop}

{\em Proof.}
For  $j\in J:=\{1,...,m \}$
we define
$$
H_{j,+}:=\{v\in V\,\,|\,\, l_j(v)=1\},
$$
$$
H_{j,-}:=\{v\in V\,\,|\,\, l_j(v)=-1\}.
$$
Let ${\cal C}$ be the set of connected components of
$V- \bigcup_{j=1}^m(H_{j,+}\cup H_{j,-})$.
For a $C\in {\cal C}$ we put 
$$
J_C:=\{j\in J \,|\, |l_j(v)|< 1\hskip 0,1cm {\rm for}\hskip 0,2cm {\rm all}
\hskip 0,2cm v\in C\}
$$
and
$$
V^C:=V\cap \bigcap_{j\in J_C}H_j.
$$
Denote by $V_C$ the complement to
$V^C$ in $V$ and let $\pi_C\, :\, V=V_C\oplus V^C\ra V_C$ be the
projection.
Since the map $V_C\ra \R^{J_C}, v\ra (l_j(v))_{j\in J_C}$,
is injective we see that $\pi_C(C)$ is a bounded open subset of $V_C$.
For $v_1\in \pi_C(C)$ we put 
$$
C(v_1):=\{v'\in V^C\,|\, v_1+v'\in C\}.
$$
The set $C(v_1)$ is a convex open subset of $V^C$. 
Let $dv_1,dv'$ be measures on $V_C$ (resp. on $V^C$), with 
$dv_1 dv'=dv$.
For all $s\in (0,1]$ we have
$$
s^{m-d}\int_{C}f(sx_0+iv)dv=s^{m-d}\int_{\pi_C(C)}\int_{C(v_1)}
f(sx_0+iv_1+iv')dv'dv_1
$$
$$
= s^{m-d^C}\int_{\frac{1}{s}\pi_C(C)}\int_{C(sv_1)}f(sx_0+isv_1+iv')dv'dv_1
$$
$$
=\int_{V}\chi_{C_s}(v)s^{m-d^C}f_{C,s}(v)dv.
$$
Here we denoted by $d^C:=\dim V^C$ and by
$$
C_s:=\{v=v_1+v'\in V_C\oplus V^C\,|\, sv_1+v'\in C\}
$$
and by $\chi_{C_s}$ the set-theoretic characteristic
function of $C_s$. We have
put for any $s\in (0,1]$ and any $v=v_1+v'\in V$
$$
f_{C,s}(v):=f(sx_0+isv_1+iv').
$$
\no
The set 
$$
{\bf K}_C:=\{ sx_0+iv_1\,|\, s\in [0,1], v_1\in \overline{\pi_C(C)}\}
$$
is contained in $T_B$ and is compact.
Further, there exist $c',c''\ge 0$ such that for all
$s\in [0,1]$ and all $v=v_1+v'\in C$ we have
$$
 \left|
\prod_{j\in J_C}(l_j(sx_0+iv)+1)\prod_{j\notin J_C}\frac{l_j(sx_0+iv)+1}{
l_j(sx_0 +iv)}\right|\le c',
$$
$$
1\le c''\cdot \left|\prod_{\j\notin J_C}\frac{1}{l_j(x_0+iv_1)}\right|.
$$
Therefore, for $s\in (0,1]$ and $v\in V$ we have
$$
|\chi_{C_s}(v)s^{m-d^C}f_{C,s}(v)|=
$$
$$
\left|\chi_{C_s}(v)g(sx_0+isv_1+iv')s^{m-d^C}
\prod_{j\in J}(1+\frac{1}{l_j(sx_0+isv_1+iv')})\right|
$$
$$
\le \chi_{C_s}(v)s^{m-d^C}\left|\prod_{j\in J_C}\frac{1}{l_j(sx_0+isv_1)}
\prod_{j\notin J_C}\frac{1}{l_j(x_0+iv_1)}\right|
c' c''\kappa ({\bf K}_C)c(v')
$$
$$
\le c'c''\kappa({\bf K}_C)\left|\prod_{j\in J}\frac{1}{l_j(x_0+iv_1)}\right|
c(v')s^{m-d^C-\#J_C}
$$
$$
\le c'c''\kappa({\bf K}_C)
\left|\prod_{j\in J}\frac{1}{l_j(x_0+iv_1)}\right|c(v'),
$$
since $m-d^C-\#J_C\ge 0$. (The constant $\kappa({\bf K}_C)$ and
the function $c\,:\,V\ra \R_{\ge 0}$ were introduced above.)

The ${\cal X}$-function corresponding to  
the cone $E^+_0\subset E_0$ (and the measure $dy$) is
given by
$$
{\cal X}_{E^+_0}(x_0+iv_1)=\prod_{j=1}^m\frac{1}{l_j(x_0+iv_1)}.
$$
Since the map from $V_C$ to $E_0$ is injective and since
$\pi_0(V_C)\cap E^+_0=\{{\bf 0}\}$ we know 
that the function
$v_1\ra {\cal X}_{E^+_0}(x_0+iv_1)$ is absolutely 
integrable over $V_C$ (by 7.2). Therefore, 
$$
v=v_1+v'\mapsto c'c''\kappa({\bf K}_C)|\prod_{j\in J}
\frac{1}{l_j(x_0+iv_1)}|c(v')
$$
is integrable over $V$. 

For a fixed $v\in V$ we consider the limit
$$
\lim_{s\ra 0}\chi_{C_s}(v)s^{m-d^C}f_{C,s}(v).
$$
The estimate above shows that this limit is $0$ if $m-d^C-\# J_C>0$.
Therefore, we assume that $m=d^C+\# J_C$. 
Then the map $V^C\ra \R^{J - J_C}$ is
an isomorphism. Since $\pi_0 (V)\cap E^+_0=\{{\bf 0}\}$, it follows
that $J_C=J$. There exists exactly one $C\in {\cal C}$ with $J_C=J$
and we denote it by $C^{\circ}$. This $C^{\circ}$ contains ${\bf 0}$ and for
all sufficiently small $s>0$ we have $s\cdot v\in C^{\circ}$, and therefore,
$v\in C^{\circ}_s$. 

Moreover, we have
$$
\lim_{s\ra 0}s^mf_{C,s}(v)=\lim_{s\ra 0}s^mg(sx_0+isv)\prod_{j=1}^m
\frac{l_j(sx_0+isv)+1}{l_j(sx_0+isv)}
$$
$$
=g({\bf 0})\prod_{j=1}^m\frac{1}{l_j(x_0+iv)}.
$$ 
Using the theorem of dominated convergence (Lebesgue's theorem),
we obtain
$$
\begin{array}{rcl}
\lim_{s\ra 0}s^{m-d}\tilde{f}_{B^+}(sx_0) & = & \lim_{s\ra 0}
\sum_{C\in {\cal C}}\frac{1}{(2 \pi )^d}\int_V \chi_{C_s}(v)
s^{m-d}f_{C,s}(v)dv\\
  &   &             \\
  & = & \frac{1}{(2\pi )^d}g({\bf 0})
\int_V \prod_{j=1}^m \frac{1}{l_j(x_0+iv)}dv\\
  &   &             \\
  & = &  g({\bf 0}){\cal X}_{\L}(\psi (x_0)).
\end{array}
$$
\hfill $\Box$

\section{Some statements on Eisenstein series}
\label{8}

\subsection*\no
{\bf 8.1}\hskip 0,5cm
Let $G$ be a semi-simple simply connected algebraic group which is
defined and split over $F$. Fix a Borel subgroup $P_0$ (defined over
$F$) and a Levi decomposition $P_0=S_0U_0$, where $S_0$ is a
maximal $F$-rational torus of $G$. Denote by 
${\bf g}$ (resp. ${\bf a}_0$)
the Lie algebra of $G$ (resp. $S_0$). 

We are going to define a certain
maximal compact subgroup ${\bf K}_G\subset G({\A})$.
This maximal compact subgroup will have the advantage that the constant term
of Eisenstein series, more precisely, certain intertwining operators, can 
be calculated explicitly, uniformly with respect to all places of $F$. 
In general, i.e., for an arbitrary maximal compact subgroup, there will be
some places where such an explicit expression is not available. In principle,
this should cause no problems. Any statement in this section should be valid
for an arbitrary maximal compact subgroup.

\subsection*\no
{\bf 8.2}
\hskip 0,5cm
Let $\Phi=\Phi(G,S_0)$ be the root system of $G$ with respect to $S_0$. 
We denote by $\Delta_0$ the basis of simple roots determined by $P_0$. 
For $\a\in \Phi$ let 
$$
{\frak g}_{\a}:= \{ X\in {\frak g}\,\,|\,\, [H,X]=\a(H)X\}
$$
be the corresponding root space.

Let $((H_{\a})_{\a\in \Delta_0}, (X_{\a})_{\a\in \Phi})$ be the 
Chevalley basis of ${\frak g}$. In  particular, this means that
$$
{\frak g}_{\a}=FX_{\a}\,(\a\in \Phi),\,\,\,\,\,\,[X_{\a},X_{-\a}]=H_{\a} \,\, 
(\a\in \Delta_0 ),
$$
$$
{\frak a}_0= \oplus_{\a\in \Delta_0}FH_{\a}.
$$
Put 
$$
{\frak g}_{\Q}=\sum_{\a\in \Delta_0}\Q H_{\a} +
\sum_{\a\in \Phi}\Q X_{\a}\subset {\frak g},
$$
This is a $\Q$-structure for
${\frak g}$ and for any $v\in \Val (F)$ the Lie algebra of $G(F_v)$ is 
${\frak g}\otimes_{\Q}F_v$.
We put 
$$
{\frak k}:=\oplus_{\a\in \Phi^+}\R(X_{\a}-X_{-\a}),
$$
$$
{\frak  p}:=\bigoplus_{\a\in \Delta_0}\R H_{\a}
\oplus\bigoplus_{\a\in \Phi^+}\R (X_{\a}+X_{-\a}),
$$
where $\Phi^+$ is the set of positive roots of $\Phi$ determined by $\Delta_0$.
Then ${\frak k}\oplus {\frak p}$ is a Cartan decomposition of 
${\frak g}_{\Q}\otimes_{\Q}\R$,  ${\frak g}_{c}:={\frak k}
\oplus i{\frak p}\subset {\frak g}_{\Q}\otimes_{\Q}\C$ is a compact form of
${\frak g}_{\Q}\otimes_{\Q}\C$ and 
${\frak g}_{\Q}\otimes_{\Q}\C={\frak g}_{c}\oplus i{\frak g}_{c} $ is a 
Cartan decomposition of ${\frak g}_{\Q}\otimes_{\Q}\C$. 

For any complex place $v$ of $F$ we define ${\bf K}_v$ to be 
$\langle \exp({\frak g}_{c})\rangle\subset G(F_v)$
(identifying $F_v$ with $\C$ via a corresponding 
embedding $F\hookrightarrow \C$).
If $v$ is a real place of $F$ we define 
${\bf K}_v=G(F_v)\cap\langle \exp({\frak g}_c)\rangle$ (identifying 
$F_v(\sqrt{-1})\simeq \C$ via a corresponding embedding $F\hookrightarrow \R$).
In this case, ${\bf K}_v$ contains $\langle \exp({\frak k})\rangle$.
In both cases ${\bf K}_v$ is a maximal compact subgroup of $G(F_v)$.

Now let $v$ be a finite place of $F$ and let ${\bf K}_v$ be the stabilizer of
the lattice 
$$
\sum_{\a\in \Delta_0}{\cal O}_v\cdot H_{\a}+
\sum_{\a\in \Phi}{\cal O}_v\cdot X_{\a}\subset {\frak g}\otimes_FF_v.
$$
By \cite{Bruhat}, sec. 3, 
Example 2, ${\bf K}_v $ is a maximal compact subgroup
of $G(F_v)$. In any case, the Iwasawa decomposition $G(F_v)=P_0(F_v){\bf K}_v$
holds (for non-archimedean $v$, cf. \cite{Bruhat}, loc. cit.).
Then ${\bf K}_G=\prod_v{\bf K}_v$ is a maximal compact subgroup of $G(\A)$ 
and $G(\A)=P_0(\A){\bf K}_G$.

\subsection*\no
{\bf 8.3}
\hskip 0,5cm
As in section 2.3 we defined for any standard parabolic subgroup
$P\subset G$
$$
H_P=H_{P,{\bf K}_G}\,:\, G({\bf A})\ra \Hom_{\C}(X^*(P)_{\C},\C)
$$
by $\langle \l, H_P(g)\rangle=\log (\prod_v|\l(p_v)|_v)$ for
$\l\in X^*(P)$ and $g=pk$,
$p=(p_v)_v\in P(\A), k\in {\bf K}_G$.

The restriction of $H_{P_0}$ to $S_0(\A)$ is a homomorphism, its kernel will
be denoted by $ S_0(\A)^1$.
The choice of a projection ${\bf G}_m(\A)\ra {\bf G}_m(\A)^1$ induces 
by means of an isomorphism $S_0\ra {\bf G}_{m,F}^{\# \Delta_0}$ a projection
$S_0(\A)\ra S_0(\A)^1$ and this in turn gives an embedding
$$
{\cal U}_0:=(S_0(\A)^1/S_0(F)(S_0(\A)\cap {\bf K}_G))^*\hookrightarrow
(S_0(\A)/S_0(F)(S_0(\A)\cap {\bf K}_G))^*.
$$
Let $(\w_{\a})_{\a\in \Delta_0}$ be the basis of $X^*(S_0)$ 
which is determined
by $\langle \w_{\a},\b\rangle =\d_{\a \b}$ for all $\a,\b\in \D_0$. 
Let $P\subset G$ be a standard parabolic subgroup. Then $\w_{\a}$ for
all $\a\in \D_P$ lifts to a character of $P$ and $(\w_{\a})_{\a\in \D_P}$
is a basis of $X^*(P)$. Put
$$
{\cal U}_P:=\{\c\,\,|\,\, \c=
\c_0\circ\prod_{\a\in \D_P}\check{\a}\circ\w_{\a} \,\,{\rm with}\,\,
\c_0\in {\cal U}_0\}.
$$
Any $\c\in {\cal U}_P$ is a character of $P(\A)/P(F)(P(\A)\cap {\bf K}_G).$
Define
$$
\phi_{\c}\,:\, G(\A)\ra S^1
$$
by $\phi_{\c}(pk)=\c(p)$ for $p\in P(\A), k\in {\bf K}_G$.
The Eisenstein series
$$
E^G_P(\l,\c,g)=\sum_{g\in P(F)\ba G(F)}\phi_{\c}(\g g)e^{\langle
\l+\rho_P,H_P(\g g)\rangle}
$$
converges absolutely and uniformly for $\Re(\l)$ contained in any compact
subset of the open cone $\rho_P+X^*(P)^+$
(cf. \cite{G}, Th\'eor\`eme III) and can be continued meromorphically to the
whole of $X^*(P)_{\C}$.

For the Eisenstein series corresponding to 
$P_0$ a proof is given in \cite{MW},
chapitre IV. In section 8.4 we will give an 
explicit expression for the Eisenstein series $E^G_P,
$ with $P\neq P_0$ as an 
iterated residue of $E^G_{P_0}$ which shows the claimed meromorphy on
$X^*(P)_{\C}$.

Let $\chi\in {\cal U}_0$. The constant term of $E_{P_0}^G (\l,\chi)$ along
$P=LU$ is by definition
$$
E^G_{P_0}(\l,\chi)_P(g)=\int_{U(F)\backslash U(\A)}E_{P_0}^G(\l,\chi,ug)du,
$$
where the Haar measure on $U(\A)$ is normalized such that $U(F)\backslash
U(\A)$ gets volume one. It is an elementary calculation to show that for any
parabolic subgroup $P\not\supsetneq P_0$ the constant term 
$E^G_{P_0}(\l,\chi)_P$ is orthogonal to all cusp forms in 
$A_0(L(F)U(\A)\ba G(\A))$ 
(cf. \cite{MW}, I.2.18, for the definition of this space). 
More precisely, for any parabolic subgroup $P\not\supseteq P_0$
the cuspidal component of $E^G_{P_0}(\l,\chi)$ along $P$ vanishes
(cf. \cite{MW}, I.3.5, for the definition of ``cuspidal component'').

By Lemme I.4.10 in \cite{MW}, the singularities of the 
Eisenstein series $E^G_{P_0}(\l,\chi)$ and the 
singularities of $E^G_{P_0}(\l,\chi)_{P_0}$
coincide. Let
$$
{\cal W}=\Norm_{G(F)}\left(S_0(F)\right)/S_0(F)
$$
be the Weyl group of $G$ with respect to $S_0$. For any $w\in {\cal W}$
we normalize the Haar measures such that 
$$
\int_{(U_0(F)\cap wU_0(F)w^{-1})\ba (U_0(\A)\cap wU_0(\A)w^{-1})}du=1
$$
and on $(U_0(\A)\cap wU_0(\A)w^{-1})\ba U_0(\A)$ we take the quotient
measure. Using Bruhat's decomposition $G(F)=\bigcup_{w\in {\cal W}}
P_0(F)w^{-1}P_0(F)$ we can calculate

\begin{equation}
\label{8.3.1}
\int_{U_0(F)\ba U_0(\A)}E^G_{P_0}(\l,\c,ug)du=\sum_{w\in {\cal W}}
c(w,\l, \c)\phi_{w\c}(g)e^{\langle w\l+\rho_0,H_{P_0}(g)\rangle},
\end{equation}

\no
where $\rho_0=\rho_{P_0}$, $(w\chi)(t)=\chi(w^{-1}tw)$ for 
all $t\in S_0(\A)$, and the functions $c(w,\l,\chi)$ are given by
$$
c(w,\l,\chi):=\int_{(U_0(\A)\cap wU_0(\A)w^{-1})\ba U_0(\A)}
\phi_{\c}(w^{-1})e^{\langle \l+\rho_P,H_{P_0}(w^{-1}u)\rangle}du
$$
(cf. \cite{MW}, Prop. II.1.7). 
They satisfy functional equations:

\begin{equation}
\label{8.3.2}
E^G_{P_0}(\l,\c,g)=c(w,\l,\c)E^G_{P_0}(w\l,w\c,g)
\end{equation}

\begin{equation}
c(w'w,\l,\c)=c(w',w\l,w\c)c(w,\l,\c)
\end{equation}

\no
(cf. \cite{MW}, Th\'eor\`eme IV.1.10).
Therefore, it suffices to calculate $c(w_{\a},\l,\c) $ for
$\a\in \D_0$ ($w_{\a}$ corresponds to the reflection along $\a$).

Put $S_{\a}=\Ker (\a)^{0}\subset S_{0}$ and $G_{\a}=Z_{G}(S_{\a})$.
The Lie algebra of $G_{\a}$ is ${\frak a}_0\oplus {\frak g}_{-\a}\oplus
{\frak g}_{\a}$. There is a homomorphism 
$\varphi_{\a}\,:\, SL_{2,F}\ra {\cal D}G_{\a}$ ($=$ derived group of $G_{\a}$)
such that $d\p_{\a}$ maps
the matrices 
$$
\left(\begin{array}{cc}0&1\\  0 & 0\end{array}\right),
\left(\begin{array}{cc}1&0\\  0 & -1\end{array}\right),
\left(\begin{array}{cc}0&0\\  1 & 0\end{array}\right)
$$
to $X_{\a}, H_{\a},X_{-\a}$, respectively. 

On $\A$ we take the measure $dx$ that is described in Tate's thesis
(then $\Vol (F\ba \A)=1$). We have
$$
c(w_{\a},\l,\chi)=\int_{\A}\phi_{\chi}(\p_{a}
( \begin{array}{cc}
\scriptstyle{0}&\scriptstyle{-1}\\
 \scriptstyle{1}&\scriptstyle{x}\end{array}) )
\exp\left(\langle \l+\rho_0,H_{P_0}(
\p_{\a}(\begin{array}{cc}
\scriptstyle{0}&\scriptstyle{-1}\\
 \scriptstyle{1}&\scriptstyle{x}\end{array}))\rangle \right)dx.
$$
It is an exercise to compute this integral. The result is

\begin{equation}
\label{8.3.3}
c(w_{\a},\l,\chi)=\frac{L(\c\circ\check{\a},\langle \l,\a\rangle)}{
L(\c\circ\check{\a},1+\langle \l,\a\rangle)}.
\end{equation}

\no
The Hecke $L$-functions are defined as follows.
Let $\c\,:\, {\bf G}_m(\A)/{\bf G}_m(F)\ra S^1$ be an unramified
character.
For any finite place $v$ we put 
$$
L_v(\c_v,s)=(1-\c_v(\pi_v)|\pi_v|^s_v)^{-1}
$$
and 
$$
L_f(\c,s)=\prod_{v\nmid \infty }L_v(\c_v,s).
$$
For any archimedean place $v$ there is a $\tau_v\in \R$ such that 
$\c_v(x_v)=|x_v|^{i\tau_v}_v$ for all $x_v\in F_v^*$. 
Then 
$$
L_v(\c_v,s):=\left\{
\begin{array}{rccl}
 \pi^{-\frac{(s+i\tau_v)}{2}}\Gamma(\frac{s+i\tau_v}{2}), & {\rm
if} & v & {\rm is}\,\,{\rm real}\\
(2\pi)^{-(s+i\tau_v)}\Gamma(s+i\tau_v), & {\rm
if} & v & {\rm is}\,\,{\rm complex}
\end{array}
\right.
$$
We define the complete Hecke $L$-function by
$$
L(\c,s)=D^{s/2}L_{\infty}(\c,s)L_f(c,s),
$$
where $L_{\infty}(\c,s)=\prod_{v|\infty}L_v(\c_v,s)$ and $D=D(F/\Q)$ is the
absolute value of the discriminant of $F/\Q$.
If the restriction of $\c$ to ${\bf G}_m({\bf A})^1$ is not trivial 
then $L(\c,s)$ is an entire function. If $\c=1$ then $ L(\c,s)$ has exactly two
poles of order one at $s=1 $ and $s=0$. To state the functional
equation we let $(\pi_v^{d_v})$  
(with $d_v \ge 0$ and $d_v=0$ for almost all $v$)
be the local discriminant of $F_v$ over 
the completion of $\Q\in F_v$ (for non-archimedean places $v$). 
Put $\delta=(\delta_v)_v\in {\bf G}_m({\bf A})$ with
$\delta_v=1$ for all archimedean places and $\delta_v=\pi_v^{d_v}$
for all non-archimedean places. Then

\begin{equation}
\label{8.3.4}
L(\c,s)=\c(\delta)L(\c^{-1},1-s).
\end{equation}

Using the functional equations (\ref{8.3.2}) and (\ref{8.3.3}) we get

\begin{equation}
\label{8.3.4-1}
c(w,\l,\chi)=\prod_{\a >0,\,\, w\l <0}\frac{L(\chi\circ\check{\a},\langle
\l,\a\rangle)}{L(\chi\circ\check{\a},1+\langle
\l,\a\rangle)}.
\end{equation}

\bigskip
\no
{\bf Proposition 8.3}\hskip 0,5cm
{\it
Let $\D_0(\chi)$ be the set of $\a\in \D_0$ such that $\chi\circ\check{\a}=1$.
Then 
$$
\prod_{\a\in \D_0(\chi)}\langle \l,\a\rangle E^G_{P_0}(\l+\rho_0,\chi)
$$
has a holomorphic continuation to the tube domain over 
$-\frac{1}{4}\rho_0 + X^*(P_0)^+$.
}
\bigskip

{\em Proof.}
For $\a\in \Phi^+- \D_0$ we have $\langle \rho_0,\a\rangle \ge 2$ and 
therefore 
$$
c(w_{\a},\l+\rho_0,\c)=\frac{L(\chi\circ\check{\a},\langle
\l+\rho_0,\a\rangle)}{L(\chi\circ\check{\a},1+\langle
\l+\rho_0,\a\rangle)}
$$
is holomorphic in this domain. If $\a\in \D_0- \D_0(\chi)$ then
$\c\circ\check{\a}$
restricted to ${\bf G}_m(\A)^1$ is nontrivial, hence
$c(w_{\a},\l+\rho_0,\c)$
is holomorphic in this domain too. 
For $\a\in \D_0(\chi)$ the function
$$
\langle \l,\a\rangle L(\chi\circ\check{\a},\langle
\l+\rho_0,\a\rangle)=\langle \l,\a\rangle L(1,1+\langle \l ,\a\rangle )
$$
is also holomorphic in this domain. This shows the holomorphy
of 
$$
\prod_{\a\in \D_0(\chi)}\langle \l,\a\rangle 
E^G_{P_0}(\l+\rho_0,\chi)_{P_0}(g)
$$
$$
= \sum_{w\in {\cal W}}\prod_{\a\in \D_0(\chi)}
\langle \l,\a\rangle c(w,\l+\rho_0,\chi)\phi_{w\chi}(g)e^{\langle
w(\l+\rho_0)+\rho_0,H_{P_0}(g)\rangle }
$$
for $\Re (\l)$ contained in  $-\frac{1}{4}\rho_0+
X^*(P_0)^+$.
By \cite{MW}, Lemme I.4.10, we conclude  that the same holds for
$\prod_{\a\in \D_0(\chi)}E^G_{P_0}(\l+\rho_0,\chi)$.
\hfill $\Box$

\subsection*\no
{\bf 8.4}
Let $P=LU$ be a standard parabolic subgroup and $\chi\in {\cal U}_P$.
For $\l\in X^*(P)_{\C}$ with $\Re (\l)$ contained in the interior
of $ X^*(P)^+ $ and $\vartheta$ contained in 
$X^*(P_0)^+ $ we have
$$
E^G_{P_0}(\vartheta+\l+\rho_0,\chi,g)
=\sum_{\g\in P(F)\ba G(F)}
\phi_{\chi}(\g g)e^{\langle \l,H_P(\g g)\rangle}\sum_{\delta\in
(L\cap P_0)(F)\ba L(F)}e^{\langle \vartheta+2\rho_0,H_{P_0}(
\delta \g g)\rangle}.
$$
Let $w_L$ be the longest element of the Weyl group of $L$ (with respect to 
$S_0$) and define
$$
c_P:=\lim_{\vartheta\ra {\bf 0},
\vartheta\in X^*(P_0)^+ } \left(\prod_{\a\in \D_0^P}
\langle \vartheta ,\a\rangle\right )c(w_L,\vartheta+\rho_0,1).
$$
By (\ref{8.3.4-1}) this limit exists and is a positive real number.

\bigskip
\no
{\bf Proposition 8.4}\hskip 0,5cm
{\it
a)  
$$
\lim_{\vartheta\ra {\bf 0}, \,\vartheta\in X^*(P_0)^+ }
 \prod_{\a\in \D_0^P}
\langle \vartheta ,\a\rangle E^G_{P_0}(\vartheta+\l+\rho_0,\chi,g)=
c_PE^G_P(\l+\rho_P,\chi,g).
$$

b) Let $P'=L'U'$ be a standard parabolic subgroup containing $P$ and 
suppose that $\chi\circ\check{\a}=1$ for all $\a\in\D^{P'}_0$.
Then $\chi\in {\cal U}_{P'}$ and for all $\l\in X^*(P')_{\C}$ we have
$$
\lim_{\vartheta\ra {\bf 0},\,
\vartheta\in X^*(P)^+}\prod_{\a\in \D_P- \D_{P'}}
\langle \vartheta ,\a\rangle E^G_{P}(\vartheta+\l+\rho_P,\chi,g)=
\frac{c_{P'}}{c_P}E^G_{P'}(\l+\rho_{P'},\chi,g).
$$
}
\bigskip 

{\em Proof.}
a) The proof rests on the fact that a 
measurable function of moderate growth on
$L(F)\ba L(\A)$ for which all cuspidal components vanish vanishes
almost everywhere. (For a proof cf. \cite{MW}, Prop. I.3.4.)
We claim that
$$
\lim_{\vartheta\ra {\bf 0},\,
\vartheta\in X^*(P_0)^+}\,\,\prod_{\a\in \D_0^P}
\langle \vartheta ,\a\rangle 
\left\{
\sum_{\delta\in (L\cap P_0)(F)\ba L(F)}
e^{\langle \vartheta+2\rho_0,H_{P_0}(
\delta \g g)\rangle}
\right\}
= c_Pe^{\langle 2\rho_P,H_{P}(\g g)\rangle }.
$$
In fact, the cuspidal components 
of both sides along all non-minimal standard
parabolic subgroups of $L$ vanish. 
To compare the constant terms along $P_0\cap
L$ we can use (\ref{8.3.1}) 
(for $L$ instead of $G$ and $P_0\cap L$ instead of
$P_0$) and the explicit expression of the functions
$c(w,\vartheta+\rho_0,1)$ to get the identity stated above (note that
$w_L\rho_0+\rho_0=2\rho_P$). 

b) Write 
$\chi=\chi_0\cdot\prod_{\a\in \D_P}\check{\a}\circ \w_{\a}$ with 
$\chi_0 \in {\cal U}_{P'}$. For $\a\in \D_P- \D_{P'}$ we have
$1=\chi\circ\check{\a}=\chi_0\circ\check{\a}$. Thus $
\chi=\chi_0\circ\prod_{\a\in \D_{P'}}\check{\a}\circ 
\w_{\a}\in {\cal U}_{P'}$.
Using a) we get
$$
\lim_{\vartheta\ra {\bf 0},\,
\vartheta\in X^*(P)^+}\prod_{\a\in \D_P-\D_{P'}}
\langle \vartheta ,\a\rangle E^G_P(\vartheta+\l+\rho_P,\chi,g)
$$
$$
= \frac{1}{c_P} \lim_{\vartheta\ra {\bf 0},
\,\vartheta\in X^*(P)^+}\prod_{\a\in \D_0^{P'}}
\langle \vartheta ,\a\rangle E^G_P(\vartheta+\l+\rho_0,\chi,g)
= \frac{c_{P'}}{c_P}E^G_{P'}(\chi+\rho_{P'},\chi,g).
$$
\hfill $\Box$

\subsection*\no
{\bf 8.5}\hskip 0,5cm
For $\c\in (S_0({\bf A})/S_0(F)(S_0({\bf A})\cap {\bf K}_G))^*$ 
and $v\in \Val_{\infty}(F) $ there is a character
$\l_v=\l_v(\c)\in X^*(S_0)_{\R}$ such that for all $x\in S_0(F_v)$
$$
\c(x)=e^{i\log(|\l_v(x)|_v)}.
$$
This gives a homomorphism 
$$
(S_0({\bf A})/S_0(F)(S_0({\bf A})\cap {\bf K}_G))^*
\ra X^*(S_0)_{\R}^{\Val_{\infty}(F)}
$$ 
$$
\l\mapsto \c_{\infty}=(\l_v(\c))_{v|\infty}
$$
which has a finite kernel. The image of ${\cal U}_0$ under this map 
is a lattice of rank 
$$
(\#\Val_{\infty}(F) -1)\dim S_0.
$$
Fix a norm $\|\cdot\|$ on $X^*(S_0)_{\R}$ and denote by the same symbol
the induced maximum norm on $X^*(S_0)_{\R}^{\Val_{\infty}(F)}$
(i.e., $\|(\l_v)_{v|\infty}\|=\max_{v|\infty}\|l_v\|$).
Let $a,b>0$ and put
$$
B_{a,b}:=\{\l\in X^*(P_0)_{\R}\,\,|\,\, -\l+\frac{a}{2}\rho_0\in
X^*(P_0)^+, \,\,\, \l+\frac{b}{2}\rho_0\in X^*(P_0)^+\}.
$$
This is a bounded convex open neighborhood of ${\bf 0}$ 
in $X^*(P_0)_{\R}$. Note that if $\l\in B_{a,b}$ then 
$w_0\l+\frac{a}{2}\rho_0
\in X^*(P_0)^+$, where $w_0$ is the longest element of ${\cal W}$.

Fix an $A>0$ such that
$$
\Re \left(\langle \l+s\rho_0,\a\rangle
(\langle \l+s\rho_0,\a\rangle-1)\right)+A\ge 1
$$ 
for all $\l\in X^*(P_0)_{\C}$ with $\Re(\l)\in B_{a,b}$, $
-1-a\le \Re(s)\le 1+b$ and $\a\in\Phi^+$.
Fix $\l\in X^*(P_0)_{\C}$ with $\Re(\l)\in B_{a,b}$ and $\c\in {\cal U}_0$.
Denote by $\Phi^+(\chi)$ the set of all positive roots $\a$ such that 
$\c\circ\check{\a}=1$. Then
$$
f_{\l,\c}(s,g):=
$$
$$
\prod_{\a\in \Phi^+(\c)}
\frac{\langle \l+s\rho_0,\a\rangle (\langle \l+s\rho_0,\a\rangle-1)}{
\langle \l+s\rho_0,\a\rangle(\langle \l+s\rho_0,\a\rangle-1)+A}
\prod_{\a >0}L_f(\chi\circ\check{\a},1+\langle \l+s\rho_0,\a\rangle)
E^G_{P_0}(\l+s\rho_0,\c,g)
$$
is for any fixed $s$ in this strip an automorphic form on 
$G(F)\ba G({\bf A})$. Indeed, we observe that all cuspidal 
components of $f_{\l,c}(s,\,\cdot\,)$ along
non-minimal standard parabolic subgroups vanish. Now we can use
(\ref{8.3.1}) and the explicit formulas for the functions
$c(w,\l,\c)$ in (\ref{8.3.4}) to see that the constant term of
$f_{\l,c}(s,\,\cdot\,)$ along $P_0$ is holomorphic for $s$ in this
domain. By Lemme I.4.10 in \cite{MW} we can conclude that the same is true
for $f_{\l,\c}(s,\,\cdot\,)$.

It is our aim to apply a version of the Phragm\'en-Lindel\"of principle
due to Rademacher (cf. \cite{Rademacher}, Theorem 2) to the function
$$
s\mapsto f_{\l,\c}(s,g)
$$
in the strip $-1-a\le \Re(s)\le 1+b$.

Using the functional equations of Eisenstein series (\ref{8.3.2}) and 
$L$-functions (\ref{8.3.4}) we get
$$
f_{\l,c}(-1-a-it,g)
=\prod_{\a\in \Phi^+(\c)}\frac{\langle \l-(1+a+it)\rho_0,\a\rangle
(\langle \l-(1+a+it)\rho_0,\a\rangle -1)}{\langle \l-(1+a+it)\rho_0,\a\rangle
(\langle \l-(1+a+it)\rho_0,\a\rangle -1)+A}                 
$$
$$
\times \prod_{\a>0}(\c\circ\check{\a})(\delta)D^{-\langle
 \l-(1+a+it)\rho_0,\a\rangle}\prod_{\a>0}L_f((w_0\chi)\circ\check{\a},
1+ \langle w_0\l+(1+a+it)\rho_0,\a\rangle)
$$
$$
\times \prod_{\a>0}
\frac{L_{\infty}\left((w_0\chi)\circ\check{\a},
1+ \langle w_0\l+(1+a+it)\rho_0,\a\rangle\right)}{
L_{\infty}((w_0\chi)^{-1}\circ\check{\a},
1-\langle w_0\l+(1+a+it)\rho_0,\a\rangle)}E^G_{P_0}(
w_0\l+(1+a+it)\rho_0,w_0\c,g).
$$
Note that $L_{\infty}(\cdots) $ is a product of $\Gamma$-functions. Using 
the functional equation of the $\Gamma$-function we can derive the following
estimate: There is $c>0$ depending only on $a$ and $b$ such that
for $\Re(\l)\in B_{a,b}$ and $\c\in {\cal U}_0$ we have
\begin{equation}
\label{8.5.1}
|f_{\l,\c}(-1-a-it,g)|
\end{equation}
$$
\le cE^G_{P_0}(\Re(w_0\l)+(1+a)\rho_0,g)\times
(1+\|\Im(\l)\|+\|\c_{\infty}\|)^{\delta_{\mu}}|1+it|^{\delta_{\mu}}.
$$
where 
$\delta_{\mu}:=\mu(2+a+b)$ and
$\mu>0$ depends only on $F$ and $G$.
Moreover, assuming $c$ to be big enough, we have also
\begin{equation}
\label{8.5.2}
|f_{\l,\c}(1+b+it,g)|\le cE^G_{P_0}(\Re(\l)+(1+b)\rho_0,g).
\end{equation}
The proof of the following lemma was suggested to us by J. Franke.

\bigskip
\no
{\bf Lemma 8.5}\hskip 0,5cm
{\it
For $\Re(\l)\in B_{a,b}, \c\in {\cal U}_0$ and $-1-a\le \s\le 1+b$
the following estimate holds:
$$
|f_{\l,\c}(\s+it,g)|
\le cE^G_{P_0}(\Re(\l)+(1+b)\rho_0,g)^{\frac{\s+1+a}{2+a+b}}
$$
$$
\times
\left\{
 E^G_{P_0}(\Re(w_0\l)+(1+a)\rho_0,g)(1+\|\Im (\l)\|+\|\chi_{\infty}\|)^{
\delta_{\mu}}|2+a+\s+it|^{\delta_{\mu}}\right\}^{\frac{1+b-\s}{2+a+b}}.
$$
}
\bigskip

{\em Proof.} This follows immediately from Theorem 2 in \cite{Rademacher}
 once we have shown that for $-1-a\le \s\le 1+b$ we have
\begin{equation}
\label{*}
|f_{\l,\c}(\s+it,g)|\le c_1e^{|t|^{c_2}}
\end{equation}
for some $c_1,c_2>0$.

By (\ref{8.5.1}) and (\ref{8.5.2}), the function
$$
s\mapsto e^{s^2}f_{\l,\c}(s,g)
$$
can be integrated over the lines 
$\Re(s)=-1-a-\epsilon$ and $\Re(s)=1+b+\epsilon$
for some $\epsilon >0$.
We claim that for all $s$ with $-1-a\le \Re(s)\le 1+b$
\begin{equation}
\label{**}
e^{s^2}f_{\l,\c}(s,g)
\end{equation}
$$
=-\frac{1}{2\pi i}\int_{\Re(z)=-1-a-\epsilon}
\frac{e^{z^2}f_{\l,\c}(z,g)}{z-s}dz +
\frac{1}{2\pi i}\int_{\Re(z)=1+b+c}
\frac{e^{z^2}f_{\l,\c}(z,g)}{z-s}dz.
$$
Denote the right-hand side by $h(s,g)$. This is a measurable function of 
moderate growth on $G(F)\ba G({\bf A})$ (cf. (\ref{8.5.1}) and (\ref{8.5.2})).
All cuspidal components of $h(s,\,\cdot\,)$ 
along non-minimal standard parabolic subgroups vanish. 
The same is true for the left-hand side.

It remains to compare the constant terms of both sides along $P_0$. By the
absolute and uniform convergence of the integrals over the vertical lines
we see that the constant term of $h(s,\,\cdot\,)$ along $P_0$  is

\begin{equation}
\label{***}
-\frac{1}{2\pi i}\int_{\Re(z)=-1-a-\epsilon}
\frac{g_{\l,\c}(z,g)}{z-s}dz+
\frac{1}{2\pi i}\int_{\Re(z)=1+b+\epsilon}
\frac{g_{\l,\c}(z,g) }{z-s}dz.
\end{equation}

where
$$
g_{\l,\c}(z,g):=
\int_{U_0(F)\ba U({\bf A})} e^{z^2} f_{\l,\c}(z,ug)du.
$$
The explicit expression of the constant term of $E^G_{P_0}$ along
$P_0$ in (\ref{8.3.1}) and uniform estimates for $L$-functions as
in \cite{Rademacher}, Theorem 5 (but for a larger strip), allow us to
conclude that (\ref{***}) is just the constant term of
$e^{s^2}f_{\l,\c}(s,\,\cdot\,)$ along $P_0$. Thus, by
Proposition I.3.4 in \cite{MW}, we have established (\ref{**}).
>From (\ref{**}) it follows that $
|e^{s^2}f_{\l,\c}(s,\,\cdot\,)|$ is bounded by some 
constant in the strip $-1-a\le \Re(s)\le 1+b$ and this in turn implies
(\ref{*}).
\hfill $\Box$.

\bigskip
\no
{\bf Proposition 8.6}\hskip 0,5cm
{\it
Let $a>0$. For any $\epsilon >0$ there exist constants
$b,c>0$ such that
for all $\l\in X^*(P_0)_{\C}$ with $\Re(\l)\in B_{a,b}$ and 
$\c\in {\cal U}_0$ we have 
$$
\left|
\prod_{\a\in \D_0(\c)}\frac{\langle 
 \l,\a\rangle}{\langle  \l,\a\rangle +1}
E^G_{P_0}(\l+\rho_0,\c)
\right|
\le c(1+\|\Im(\l)\|+\|\chi_{\infty}\|)^{\epsilon}.
$$
}
\bigskip

{\em Proof.}
Note that if $\a\in \Phi^+ - \D_0$ then $\langle \rho_0,\a\rangle\ge 2$
and hence
$\langle \l+\rho_0,\a\rangle (\langle \l+\rho_0,\a\rangle -1)$ does not
vanish for $\Re(\l)\in B_{a,b}$ and $b>0$ sufficiently small. 
For such $b$ there is a constant $c_1$ such that
$$
\left|\prod_{\a\in \D_0(\c)}\frac{\langle 
 \l,\a\rangle}{\langle  \l,\a\rangle +1}E^G_{P_0}(\l+\rho_0,\c)\right|\le
c\left|f_{\l,\c}(1,1_G)\right|.
$$
Now we use the estimate for $|f_{\l,\c}(1,1_G)|$ in Lemma 8.5 and
require that $\mu b\le \epsilon$. This gives the desired result.
\hfill $\Box$

\bigskip\no
{\bf Proposition 8.7}\hskip 0,5cm
{\it
Let $P$ be a standard parabolic subgroup of $G$. Let $a,\epsilon >0$.
Then there exist $b,c >0$ such that for all $\c\in {\cal U}_P$ and 
$\l\in X^*(P)_{\C}$ with $-\Re(\l)+\frac{a}{2}\rho_0\in X^*(P)^+$ and
$\Re(\l)+\frac{b}{2}\rho_0\in X^*(P)^+$ we have
$$
\left|
\prod_{\a\in \D_P(\c)}\frac{\langle 
 \l,\a\rangle}{\langle  \l,\a\rangle +1}
E^G_{P}(\l+\rho_P,\c)
\right|
\le c(1+\|\Im(\l)\|+\|\chi_{\infty}\|)^{\epsilon},
$$
where $\D_P(\c)=\D_0(\c)\cap \D_P$.
}
\bigskip

{\em Proof.}
By the preceding proposition, there exist $b,c'>0$ such that
for all 
$\l'\in X^*(P_0)_{\C}$ with $\Re(\l')\in B_{a,b}$
we have
$$
\left|\prod_{\a\in \D_0(\c)}\frac{\langle 
 \l',\a\rangle}{\langle  \l',\a\rangle +1}E^G_{P_0}(\l'+\rho_0,\c)\right|
\le c'(1+\|\Im(\l')\|+\|\chi_{\infty}\|)^{\epsilon}.
$$
Note that $\chi\circ \check{\a}=1 $ for all $\a\in \D_0^P$ and hence
$\D_0(\c)=\D_P(\c)\cup \D_0^P$.
Now let $\l\in X^*(P)_{\C}$ be as in the proposition, i.e.,
$-\Re(\l)+\frac{a}{2}\rho_P\in X^*(P)^+$ and 
$\Re(\l)+\frac{b}{2}\rho_0\in
X^*(P)^+$. 
Then for all sufficiently small $\vartheta\in X^*(P_0)^+$ we have
$$
-(\vartheta +\Re(\l))+\frac{a}{2}\rho_0\in X^*(P_0)^+
$$
$$
\vartheta +\Re(\l)+\frac{b}{2}\rho_0\in X^*(P_0)^+,
$$
i.e.,
$\vartheta +\Re(\l)\in B_{a,b}$.
Hence for those $\vartheta$ we have
$$
\left|
\prod_{\a\in \D_0(\c)}\frac{\langle \vartheta +\l,
\a\rangle}{\langle \vartheta+ \l,\a\rangle +1}
E^G_{P_0}(\vartheta+\l+\rho_0,\c)
\right|
\le c'(1+\|\Im(\l)\|+\|\chi_{\infty}\|)^{\epsilon}.
$$
Letting $\vartheta $ tend to ${\bf 0}$ and using Proposition 
8.4 we can conclude that 
$$
\left|
\prod_{\a\in \D_P(\c)}\frac{\langle \l,
\a\rangle}{\langle \l,\a\rangle +1}E^G_{P}(\l+\rho_P,\c)
\right|
\le \frac{c'}{c_P}(1+\|\Im(\l)\|+\|\chi_{\infty}\|)^{\epsilon}.
$$
\hfill $\Box$

\end{document}